\documentclass[12pt]{amsart}

\usepackage{amsmath,amsthm,amsfonts,amscd,amssymb,pb-diagram,epsfig,graphicx,verbatim}
\usepackage[all]{xy}

\theoremstyle{plain}

\newtheorem{theorem}{Theorem}[subsection]
\newtheorem*{theorem*}{Theorem}
\newtheorem*{corollary*}{Corollary}
\newtheorem{corollary}[theorem]{Corollary}

\newtheorem*{lemma*}{Lemma}
\newtheorem{lemma}[theorem]{Lemma}

\newtheorem*{claim*}{Claim}
\newtheorem*{proof*}{proof}
\newtheorem*{prop*}{Proposition}
\newtheorem{proposition}[theorem]{Proposition}

\newtheorem{defn/thm}[theorem]{Definition/Theorem}
\newtheorem{defn/prop}[theorem]{Definition/Proposition}
\newtheorem{definition/lemma}[theorem]{Definition/Lemma}
\newtheorem{example/lemma}[theorem]{Example/Lemma}

\theoremstyle{definition}
\newtheorem{remark}[theorem]{Remark}
\newtheorem*{example*}{Example}
\newtheorem{example}[theorem]{Example}
\newtheorem{defn}[theorem]{Definition}
\newtheorem*{question*}{Question}

\newcommand{\cS}{\mathcal{S}}
\newcommand{\cI}{\mathcal{I}}
\newcommand{\cL}{\mathcal{L}}

\newcommand{\cO}{\mathcal{O}}

\newcommand{\M}{\overline{M}}

\newcommand{\ra}{\rightarrow}
\newcommand{\hra}{\hookrightarrow}

\newcommand{\ov}{\overline}
\newcommand{\til}{\widetilde}

\newcommand{\mc}[1]{\mathcal #1}

\newcommand{\Z}{\mathbb{Z}}
\newcommand{\Q}{\mathbb{Q}}
\newcommand{\bP}{\mathbb{P}}
\newcommand{\A}{\mathbb{A}}
\newcommand{\C}{\mathbb{C}}

\newcommand{\bd}{\begin{equation*} \begin{diagram}}
\newcommand{\ed}{\end{diagram} \end{equation*}}
\newcommand{\no}{\node}
\newcommand{\arr}{\arrow}

\begin{document}

\author{L. Chen, A. Gibney, and D. Krashen}

\title{Pointed trees of projective spaces \thanks{The authors were
supported during this work by National Science Foundation under
agreements DMS-0432701, DMS-0509319, and DMS-0111298 respectively.} }

\maketitle

\begin{abstract}
We introduce a smooth projective variety $T_{d,n}$ which compactifies
the space of configurations of $n$ distinct points on affine $d$-space
modulo translation and homothety. The points in the boundary correspond
to $n$-pointed stable rooted trees of $d$-dimensional projective spaces,
which for $d = 1$, are $(n+1)$-pointed stable rational curves.  In
particular, $T_{1,n}$ is isomorphic to $\overline{M}_{0,n+1}$, the
moduli space of such curves.  The variety $T_{d,n}$ shares many
properties with $\overline{M}_{0,n}$. For example, as we prove, the
boundary is a smooth normal crossings divisor whose components are
products of $T_{d,i}$ for $i < n$, it has an inductive construction
analogous to but differing from Keel's for $\overline{M}_{0,n}$ which
can be used to describe its Chow groups, Chow motive and Poincar\'e
polynomials, generalizing \cite{Keel,Man:GF}. We give a presentation of
the Chow rings of $T_{d,n}$, exhibit explicit dual bases for the
dimension $1$ and codimension $1$ cycles.  The variety $T_{d,n}$ is
embedded in the Fulton-MacPherson spaces $X[n]$ for \textit{any} smooth
variety $X$ and we use this connection in a number of ways.  For
example, to give a family of ample divisors on $T_{d,n}$ and to give an
inductive presentation of the Chow groups and the Chow motive of $X[n]$
analogous to Keel's presentation for $\overline{M}_{0,n}$, solving a
problem posed by Fulton and MacPherson.

\end{abstract}

\section{Introduction}
\label{intro}
Fix an arbitrary ground field $k$. By a variety over $k$ we mean a
reduced (but not necessarily integral), seperated scheme of finite type over
$k$.

Let $TH_{d,n}$ denote the space of configurations of $n$ distinct points
on affine $d$-space up to translation and homothety. Equivalently, this
may be regarded as the space of embeddings of a hyperplane and $n$
distinct points not lying on the hyperplane in projective $d$-space, up
to projective automorphisms. When $d=1$, this is the moduli space
$M_{0,n+1}$ of pointed rational curves. In this paper, we introduce and
study varieties $T_{d,n}$ which compactify $TH_{d,n}$ for $d\geq 1,n\geq
2$.  We prove that $T_{d,n}$ is a smooth, projective, irreducible,
rational variety of dimension $dn-d-1$ (c.f.  Corollary
\ref{blowupcor1}). The points of $T_{d,n}$ are in one to one
correspondence with \emph{stable n-pointed rooted trees} of
$d$-dimensional projective spaces (Definition \ref{tree_def}, Theorem
\ref{tree_thm}).  These pointed trees of projective spaces are
higher-dimensional analogs of stable pointed rational curves.  Indeed,
$T_{1,n} \cong \M_{0,n+1}$ (Theorem \ref{t1n}).  Remarkably, $T_{d,n}$
seems to share nearly all of the combinatorial and structural advantages
of $\M_{0,n}$.

There has been much interest in possible higher-dimensional
generalizations of $\M_{0,n}$.  For example, Kapranov's Chow quotients
compactify the moduli spaces of ordered $n$-tuples of hyperplanes in
$\bP^{r}$ in general linear position in \cite{Kap:CQ}. These are
isomorphic to $\M_{0,n}$ when $r=1$.  Hacking, Keel, and Tevelev defined
and studied another compactification of hyperplane arrangements in
projective space as the closure inside the moduli spaces of stable pairs
(see \cite{KT:CQ,HKT}).  These spaces are reducible, with Kapranov's
Chow quotients as one component.

The Fulton-MacPherson space $X[n]$, defined as a compactification of the
space of configurations of $n$ distinct points on a smooth variety $X$,
is also a kind of higher-dimensional analog of $\M_{0,n}$.  In
particular, $\bP^1[n]$ is birational, although not isomorphic, to
$\M_{0,n+3}$ (see \cite{FM}).

From a different perspective, one may show that $\M_{0,n+1}$ may be
viewed as a subscheme of $X[n]$ for any smooth curve $X$. In more
generality, $T_{d,n}$ arises as a subscheme of $X[n]$ for {\textit{any}}
smooth variety $X$ of dimension $d$.

Moduli spaces of pointed rational curves are intimately related to
moduli spaces of curves of higher genus $g$.  Namely, by attaching
curves in various ways, one can define maps from $\M_{0,n+1}$ to the
boundary of $\M_{g,m}$.  It has been possible to reduce important
questions about the birational geometry of $\M_{g,n}$ for $g > 0$ to
moduli of pointed rational curves \cite{GKM}.  Analogously, by using
various attaching maps, the variety $T_{d,n}$ maps to the boundary of
the Fulton-MacPherson space $X[n]$. In fact, $T_{d,n}$ is a fiber of the
natural projection map from the boundary component $D(N) \subset X[n]$
to $X$ (Definition \ref{def_tdn}).  We exploit this fundamental fact to
study $T_{d,n}$ as well as to answer a question posed by Fulton and
MacPherson about the Chow groups of $X[n]$ (cf. Theorem
\ref{fm_chow_groups}).

Another generalization of the moduli space of curves is the stack
$\M_{g,n}(X,\beta)$ of stable maps from $n$-pointed stable curves of
genus $g$ to a variety $X$.  When $X$ is a point and $g=0$, we recover
$\M_{0,n}$.  Stable maps have been particularly studied because their
Chow rings determine the Gromov-Witten invariants of the variety $X$.
The work of Oprea and Pandharipande, for example, show that the
combinatorial structure of $\M_{0,n}$ plays a major role in the
understanding the intersection theory of the more general spaces
(\cite{Opr:Div,Opr:TautCl,Opr:TautSM,Pand}). It is conceivable that the
$T_{d,n}$ could be used to study moduli of stable maps from higher
dimensional varieties.

The authors would like to thank P. Deligne and R. Pandharipande for
helpful conversations.

\subsection{Summary of Results}

\ \ \\

\noindent
\textit{Inductive construction, Chow groups, Chow motives, Poincar\'e
polynomials, and functor of points}

We describe a functor which represents $T_{d,n}$ in Section
\ref{construction} (cf. Proposition \ref{sheaf_functor} and Lemma
\ref{triv_tdn}) and use it to prove that $T_{d,n}$ can be constructed
inductively.  More specifically, we prove:

\begin{theorem*}[\ref{main_blowup}]
The variety $T_{d,n}$ may be constructed as the result of a sequence of
blowups of a projective bundle over $T_{d,n-1}$.
\end{theorem*}
In particular, this gives a construction of $\M_{0,n}$ which differs
from previous constructions of Keel and of Kapranov
(\cite{Keel,Kap:M0n}).  We use this construction to obtain an inductive
presentation the Chow groups and the Chow motive of $T_{d,n}$ (Section
\ref{applications}) and a description of its Poincar\'e polynomial
(Section \ref{poincare}).

\bigskip

\noindent
\textit{Ample divisors}

Using the embedding of $T_{d,n}$ as a closed subvariety of the
Fulton-MacPherson configuration space $X[n]$ for a smooth
$d$-dimensional variety $X$, we exhibit a family of ample divisors for
$T_{d,n}$ in Theorem \ref{ample}.

\bigskip

\noindent
\textit{Boundary of $T_{d,n}$, and stratification}

We characterize the boundary of the $T_{d,n}$, showing it is composed of
smooth normal crossings divisors which are (isomorphic to) products of
smaller $T_{d,i}$.  In particular, for each $S \subsetneq N$, there is a
nonsingular divisor $T_{d,n}(S) \subset T_{d,n}$ such that the union of
these divisors $T_{d,n}(S)$ forms the boundary $T_{d,n} \setminus
TH_{d,n}$.  Any set of these divisors meets transversally.  An
intersection of divisors $$T_{d,n}(S_1) \cap \cdots \cap T_{d,n}(S_r)$$
is nonempty exactly when the sets $S_i$ are \emph{nested }; each pair is
either disjoint, or one is contained in the other.  Moreover, the
boundary components $T_{d,n}(S)$ are products.  Namely, $$T_{d,n}(S)
\cong T_{d,n-|S|+1} \times T_{d,|S|},$$ for $S \subsetneq N$, $|S| > 1$
(Theorem \ref{main_blowup}, part \ref{products}).

More generally, as in the case $d=1$, the $T_{d,n}$ are stratified by
the (closure) of the locus of points corresponding to varieties having
$k$ distinct components.  There is a natural divisor class $\delta_N$
on $T_{d,n}$; for $d=1$, $\delta_N=-\psi_{n+1}$ (beginning of Section
\ref{pairing_section}). We give a simple presentation for the Chow
ring of $T_{d,n}$ in terms of the $\delta_N$ and the boundary classes:
$$A^*(T_{d,n}) \cong \mathbb{Z}[\{\delta_S\}_{S\subset N, 2 \le
|S|}]/I_{d,n},$$ where the ideal $I_{d,n}$ is generated by two simple
types of relations (Theorem \ref{Chowring}).  As in the case of
$T_{1,n}\cong\M_{0,n+1}$, there are natural maps between the spaces
given by dropping points.  That is, for every $i \in \{1, \ldots,
n\}$, there is a natural projection $T_{d,n} \ra T_{d,n-1}$ given by
``dropping the $i$th point,'' (Remark \ref{point_drop}).  We prove
that $T_{d,n}$ is an HI space (Corollary \ref{tdn_HI}) when defined
over $\C$. That is, $H^{2*}(T_{d,n})\cong A^*(T_{d,n}).$

We describe a relationship between boundary divisors and certain
explicitly given one-cycles which yields an integer pairing between
divisors and curves on $T_{d,n}$ (Theorem \ref{pairing}).  These
classes form a basis for 1-cycles modulo rational equivalence
(Corollary \ref{intbasis}).  We also give an explicit conjectural pairing
between cycles of complementary dimension on $T_{d,n}$ (cf. Section
\ref{conjectural_pairing}).

\section{The closed points of $T_{d,n}$} \label{closed_points}

In this section, we give a geometric description of the closed points of
$T_{d,n}$ in terms of isomorphism classes of $n$-pointed rooted trees of
$d$-dimensional projective spaces. 

Choose a pair of smooth $d$-dimensional varieties $X_1, X_2$, a point $p
\in X_1$ and a subvariety $H \subset X_2$ such that $H \cong \bP^{d-1}$.
Let $Y$ be the blowup of $X_1$ at $p$. We may form a new variety, which
we will denote by $X_1 \#_{p,H} X_2$ by identifying the exceptional
divisor in $Y$ with the subvariety $H \in X_2$. In the case $d = 1$ this
correponds to attaching two curves together by identifying a point on
one with a point on the other.

To describe a tree of $d$-dimensional projective spaces, we will use 
trees as book keeping devices.  Recall that a rooted tree is a graph
without cycles and with a distinguished vertex. We will use the notation
$G = (V_G, E_G, v_G)$ where $V_G$ is a set of verticies, $v_G \in V_G$
is a distinguished vertex called the root and $E_G \subset V_G \times
V_G$ is the set of egdes to denote such an object. Recall that given a
rooted tree $G$, there is a natural partial order on $V_G$ in which the
root $v_G$ is the initial or smallest element. Given $w < w'$ we
say that $w'$ is a descendant of $w$. In the case that $w < w'$ and
there is no vertex $w''$ with $w < w'' < w'$, we say that $w'$ is a
daughter of $w$ and that $w$ is the parent of $w'$.

We define $d$-dimensional gluing data for a tree $G$ to be a collection
of projective spaces $X_w \cong \bP^d$ for each $w \in V_G$ together
with a rule which associates to each vertex $w \in V_G$ a hyperplane
$H_w \subset X_{w}$ and to each pair $w, w'\in V_G$, where $w'$ is a
daughter of $w$, a point $p(w, w') \in X_w$ such that the points $p(w,
w') \in X_w$ are all distinct as $w'$ varies over the daughters of $w$,
and do not lie on the hyperplane $H(w)$. We denote this data by $(X, p,
H)$. Given such gluing data, we may define a variety $\underset{p,
H}{\#}X$ inductively on the order of $V$ as follows:
\begin{enumerate}
\item
If $|V| = 1$, then $\underset{p, H}{\#}X = X_v = \bP^d$,
\item
If $|V| = n+1$, choose a vertex $w \in V_G$ with no daughters, and let
$w' \in V_G$ be the parent of $w$. Let $G'$ be the tree obtained by
removing $w$ and all edges incident with $w$.  Let $X', p', H'$ be the
restrictions of the functions $X, p, H$ to $G'$.  Then 
\[\underset{p, H}{\#}X = \left(\underset{p', H'}{\#}X'\right)
\#_{p(w', w), H_w} X_w.\]
\end{enumerate}

We note that in the variety $\til{X} = \underset{p,H}{\#}X$, each
component is a blowup of one of the varieties $X_w$ and consequently
there is a 1-1 correspondence between the components of $\til{X}$ and
the verticies $V_G$. The singular locus of $\til{X}$ is exactly the
intersections of the different components. Each component has a
distinguished hyperplane $H_w$, which is in the singular locus of
$\til{X}$ unless $w = v_G$ is the root of $G$.

\begin{defn}
A {\emph{rooted tree of $d$-dimensional projective spaces}} (a
\emph{$d$-RTPS}) is a connected variety $Z$ together with a closed
embedding $f : \bP^{d-1}\hra Z$ (called the \emph{root}) such that
there is some rooted tree $G$ and gluing data $(X, p, H)$ such that $Z
\cong \underset{p, H}{\#}X$, and $f$ defines an isomorphism of
$\bP^{d-1}$ with $H_{v_G}$. 
\end{defn}

\begin{defn} \label{tree_def}
An \emph{$n$-pointed $d$-RTPS}  $(\bP^{d-1} \hra Z, p_1, \ldots, p_n)$ 
is a $d$-RTPS with distinct \emph{marked points} $p_1,\ldots,p_n\in Z$
such that:
\begin{enumerate}
\item
$p_i$ is not in the singular locus of $Z$,
\item
For all $i$, $p_i$ does not lie in (the image of) the root.
\end{enumerate}
\end{defn}

\begin{defn}
$(\bP^{d-1} \hra Z, p_1, \ldots, p_n)$  is \emph{stable} if each
component $W \subset Z$  contains at least two distinct \emph{markings},
where a marking is either a marked point $p_i$ or an exceptional
divisor. Note that each exceptional divisor corresponds to a daughter of
the vertex corresponding to $W$.
\end{defn}

Note that this agrees with the situation for a stable pointed rational
curve. Although the standard definition in this case is requires $3$
markings, in our general situation we do not count the hyperplanes $H_w$
as markings. Since each component has exactly one such hyperplane, this
shows that our definition is specializes to the standard one.

\begin{defn}
Two $n$-pointed rooted trees of $d$-dimensional projective
spaces $(\bP^{d-1} \hra Z, p_1, \ldots, p_n)$ and $(\bP^{d-1} \hra Z',
q_1, \ldots, q_n)$ are \emph{isomorphic} if there is an isomorphism of
algebraic varieties $f: Z \ra Z'$ such that $f(p_i) = q_i$ and the
following diagram commutes:
\bd
\no[2]{\bP^{d-1}} \arr{sw} \arr{se}\\
\no{Z} \arr[2]{e} \no[2]{Z'.}
\ed
\end{defn}

The following proposition is easy to verify explicitly:
\begin{proposition}
An $n$-pointed rooted tree of $d$-dimensional projective spaces is
stable if and only if it has no nontrivial automorphisms.
\end{proposition}

As a simple example, consider a stable $(n+1)$-pointed rational curve
consisting of three components $Z_i \cong \bP^1$ such that $Z_2$ and $Z_3$ are
attached to $Z_1$ by identifying for each $i \in \{2,3\}$ a point $h_i \in
X_i$ with $e_i \in X_1$.  Suppose there are $s \ge 2$ points $p_1$,$\ldots,p_s
\in X_2 \setminus \{h_2\}$ and $n-s \ge 2$ points $p_{s+1}$,$\ldots p_n \in
X_3 \setminus \{h_3\}$ and the $(n+1)$-st point $p_{n+1}$ is on $Z_3 \setminus
\{e_2, e_3\}$. The curve is a tree of projective lines and is illustrated
below in figure \ref{stable_lines}. We call $p_{n+1} \in X_1$ the root of the
tree.

More generally, let $Z_1=Bl_{\{q,q'\}}\bP^d$ be the blow up of
$\bP^d$ at the points $q$ and $q'$ with exceptional divisors $E_{2}$
and $E_{3}$, and let $H_1 = \bP^{d-1} \hra Z_1\setminus\{q,q'\}$ be an
embedded hyperplane.  Let $Z_2$ be isomorphic to $\bP^d$ with fixed
marked points $p_1$,$\ldots,p_s \in Z_2$ and a fixed embedded hyperplane $H_2
= \bP^{d-1} \hra Z_2 \setminus \{p_1,\ldots, p_s\}$.  Finally, let
$Z_3$ be isomorphic to $\bP^d$ with marked points
$p_{s+1}$,$\ldots$,$p_n \in Z_3$ and an embedded hyperplane $H_3 =
\bP^{d-1} \hra Z_3 \setminus \{p_{s+1},\ldots, p_n\}$. Let $H_2$ be
identified pointwise with $E_2$ and $H_3$ with $E_3$.  When the components are
attached, they from a tree.  We call the embedded hyperplane $H_1 =
\bP^{d-1} \subset X_1$ the root of the tree, shown below in figure
\ref{stable_planes} when $d = 2$.  If $s \ge 2$ and $n-s \ge 2$, the tree is
stable; it has no nontrivial automorphisms fixing the embedded hyperplanes
pointwise which preserve the marked points.

\medskip

\begin{figure}[h]
  \hfill
  \begin{minipage}[t]{.45\textwidth}
    \begin{flushleft}  
      \epsfig{file=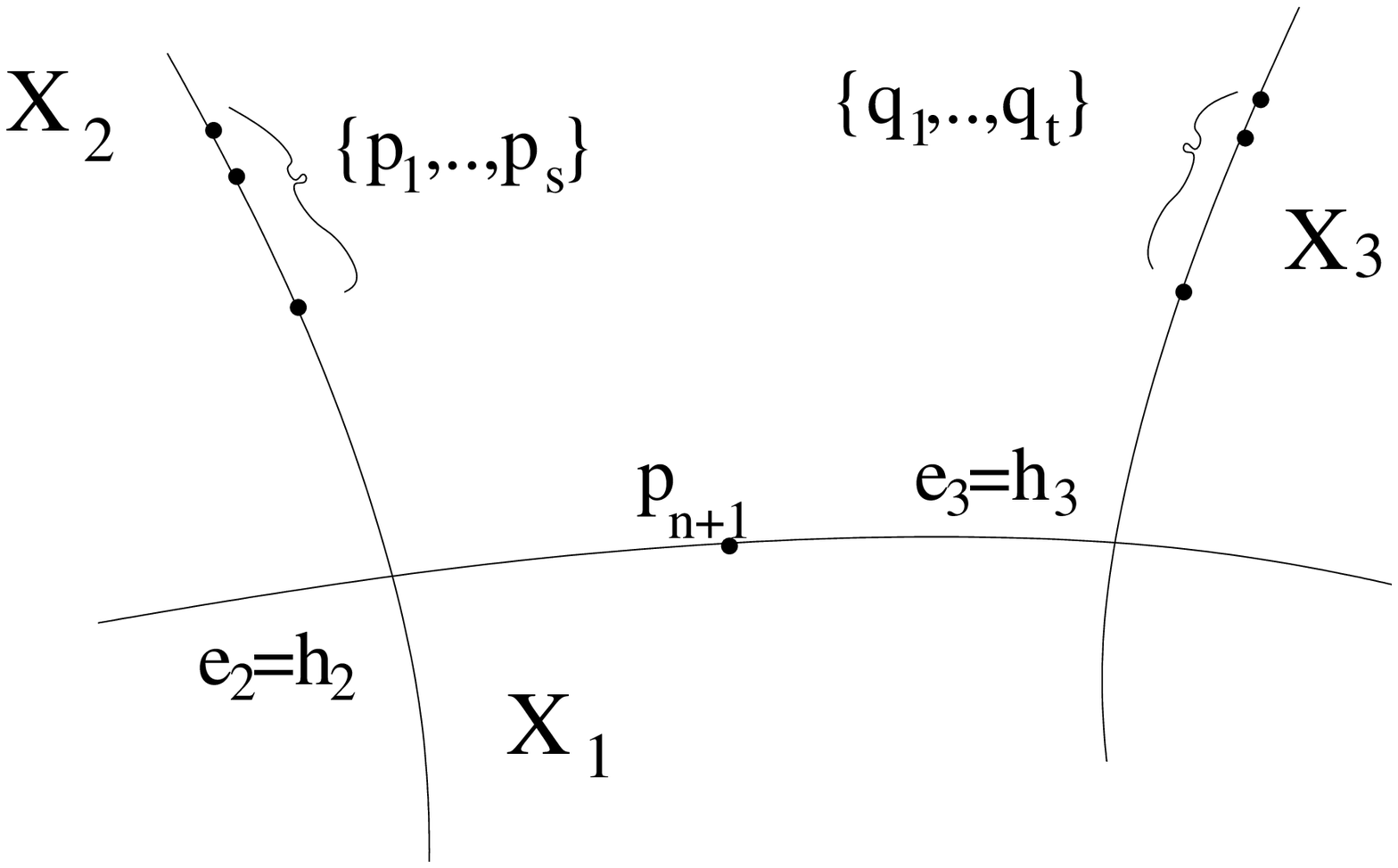, scale=0.24}
      \caption{d = 1}
      \label{stable_lines}
    \end{flushleft}
  \end{minipage}
  \hfill
  \begin{minipage}[t]{.45\textwidth}
    \begin{center}  
      \epsfig{file=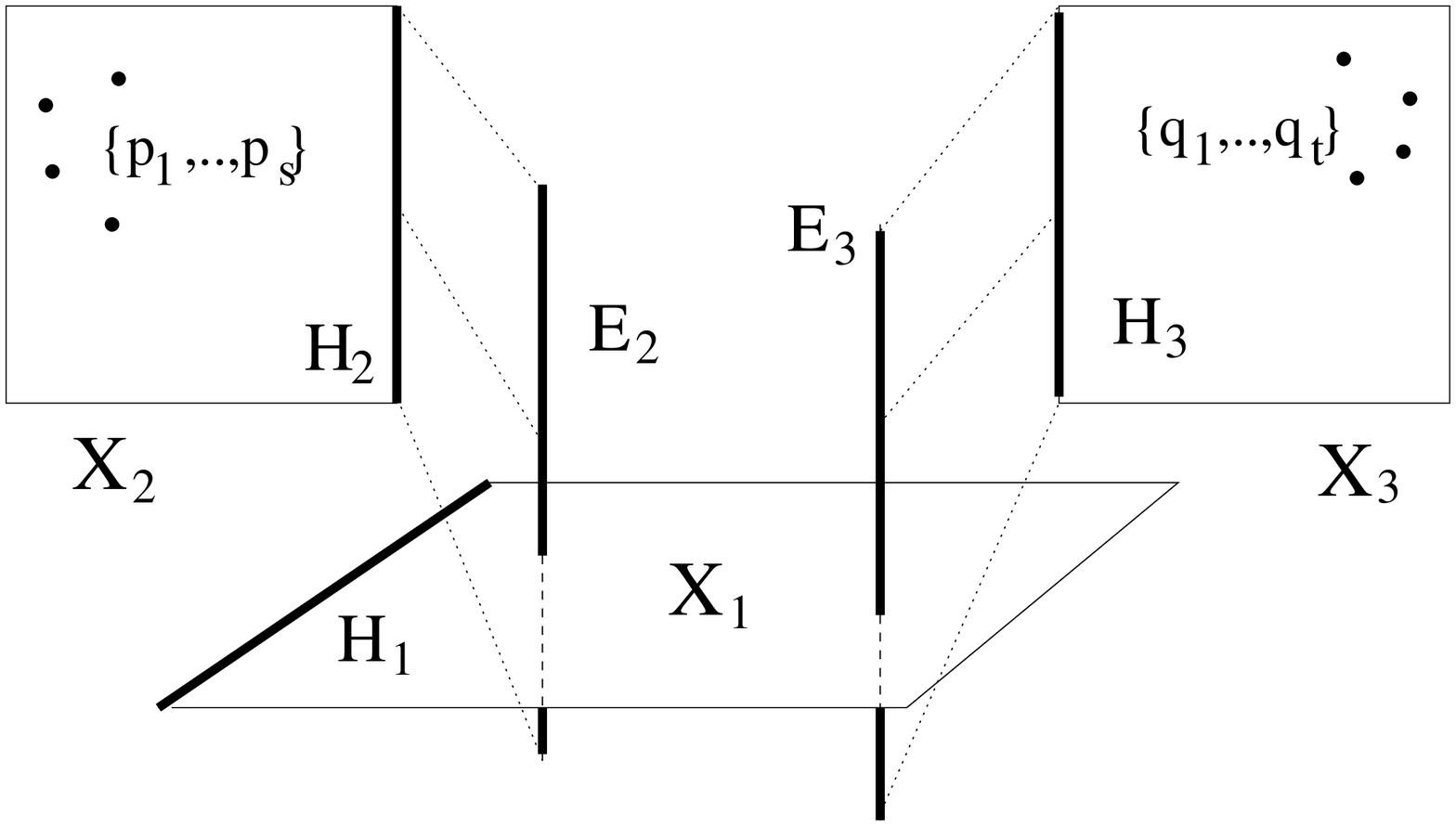, scale=0.24}
      \caption{d = 2}
      \label{stable_planes}
    \end{center}
  \end{minipage}
  \hfill
\end{figure}

\section{Definition and Inductive Construction of $T_{d,n}$}
\label{construction}

In this section we define the variety $T_{d,n}$ as an abstract variety.
This is done by using the construction of the Fulton-MacPherson
configuration space as in \cite{FM}. 

\subsection{Inductive Construction of the Fulton-Macpherson Space}

Let $X$ be a smooth variety of dimension $d$, and let $x \in X(k)$. As
in \cite{FM}, we let $X[n]$ denote the Fulton-MacPherson configuration
space of $n$ points on $X$, whose construction we recall below (with
some minor changes in notation). Set $N = \{1, \ldots, n\}$. The space
$X[n]$ comes with a morphism $X[n] \ra X^n$. For every subset $S \subset
N$ with $|S| \geq 2$, Fulton and MacPherson define a codimension $1$
smooth subvariety $D(S) \subset X[n]$ which maps into the diagonal
$\Delta_S = \{(x_i) \in X^n | x_i = x_j \text{ for } i,j \in S\}$. In
particular, we have a morphism $\pi : D(N) \ra X \cong \Delta_N\subset
X^n$.

\begin{defn} \label{def_tdn}
$T_{d,n}^{X,x} = \pi^{-1}(x)$.
\end{defn}

We shall prove that this definition does not depend on the smooth
variety $X$ or on the point $x \in X(k)$. Thus we simply write $T_{d,n}$
for $T_{d,n}^{X,x}$.  To show this, we describe the functor which it
represents  later in this section and show that this functor is
independent of our choices (see Definition \ref{sheaf_functor}). We
also show that the points of $T_{d,n}$ correspond to the $n$-pointed
stable $d$-RTPS's from the previous section (Theorem \ref{tree_thm}).

In order to set notation and motivate our work on the spaces $T_{d,n}$,
we now recall Fulton and MacPherson's construction of $X[n]$. 

The construction of these spaces is given inductively.  It will be
notationally convenient to let $N = \{1, \ldots, n\}$, and we may
occasionally write
$X[N]$ to mean $X[n]$. For a subset $S \subset N$ we let $S^+$ be the
subset $S \cup \{n+1\} \subset \{1, \ldots, n+1\}$. In particular, $N^+
= \{1, \ldots, n+1\}$. 

At the $n$'th step in the process, we will have constructed:
\begin{enumerate}
\item a space $X[n]$,
\item a morphism $\pi_n : X[n] \ra X^n$ (which we will write as $\pi$
when $n$ is understood),
\item for each subset $S \subset N$ with $|S| \geq 2$, a divisor $D(S)
\subset X[n]$.
\end{enumerate}

We begin by giving the definition of the first two spaces directly.  For
$n = 1$, we set $X[1]=X$. For $n = 2$, we let $X[2]=Bl_{\Delta}(X\times
X)$ be the blowup of $X\times X$ along the diagonal $\Delta$. We define
$D(\{1,2\})$ to be the exceptional divisor of this blowup, and let $\pi
: X[2] \ra X^2$ be the blowup map. 

To go from $X[n]$ to $X[n+1]$ requires a series of steps in itself.
We'll construct a sequence of smooth varieties:
\begin{multline*}
X[n,n] \stackrel{\rho_n}{\longrightarrow}
X[n,n-1] \stackrel{\rho_{n-1}}{\longrightarrow} \cdots \longrightarrow
X[n,k+1] \stackrel{\rho_{k+1}}{\longrightarrow} X[n,k] \longrightarrow
\\ \cdots \longrightarrow X[n,1] \stackrel{\rho_1}{\longrightarrow}
X[n,0],
\end{multline*}
so that $X[n,n] = X[n+1]$. We will define these varieties $X[n,k]$
inductively with respect to $k$. At each step, we will construct for $0
\leq k \leq n$:

\begin{enumerate}
\item a smooth variety $X[n,k]$,
\item a morphism $\rho_k : X[n,k] \ra X[n,k-1]$ when $k > 0$,
\item smooth subvarieties $X[n,k](S')$ for each subset $S' \subset N^+$
with at least two elements.
\end{enumerate}

In the case $k = 0$, we set $X[n,0] = X[n] \times X$. For $S' \subset N
\subset N^+$, we define $X[n,0](S') = D(S') \times X$. Let $p_i : X[n]
\ra X^n \ra X$ be the composition of $\pi_n$ with the $i$'th projection
map. We define $X[n,0](\{i\}^+)$ to be the graph of $p_i$ - in other
words, it is the image of the morphism $id \times p_i : X[n] \ra X[n]
\times X$. Now suppose that $S \subset N$, $|S| \geq 2$. It turns out
that if $i, j \in S$, then if we let $\Gamma_i= id \times p_i
:D(S)\longrightarrow D(S)\times X$, the images $\Gamma_i(D(S))$ and
$\Gamma_j(D(S))$ are isomorphic. We denote these common maps by
$\Gamma_S$ and define $X[n,0](S^+)$ to be the image $\Gamma_S(D(S))$.

The variety $X[n,1]$ is defined to be the blowup of $X[n,0]$ along the
subvariety $X[n,0](N^+)$. For $S' \neq N^+$, the variety $X[n,1](S')$ is
defined to be the proper transform of $X[n,0](S')$, and we define
$X[n,1](N^+)$ to be the exceptional divisor.  We have the following
pullback diagram:
\[\xymatrix{
X[n,1](N^+) \ar@{=}[r] & P(\mathcal{N}_{N}) \ar[d] \ar[r] 
& X[n,1] \ar[d]_{\rho_1} \\                                
& X[n,0](N^+) \ar[r] & X[n,0] \ar@{=}[r] & X[n] \times X,
}\]

where $\mathcal{N}_N=N_{X[n,0](N^+)}X[n,0]$ is the normal bundle of
$X[n,0](N^+)$ in $X[n,0]$ and $X[n,1](N^+) = P(\mathcal{N}_N)$ is the
exceptional divisor of the blowup. It will be useful to consider the
first Chern class of this bundle, and so we will set $l_N =
c_1(\mathcal{O}_{\mathcal{N}_N}(1))$.

Once $X[n,k]$ has been constructed for $k \geq 1$ along with its
subschemes $X[n,k](S')$, the variety $X[n,k+1]$, together with its
morphism $\rho_{k+1}:X[n,k+1] \rightarrow X[n,k]$, is defined to be the
blowup along the disjoint union of the subvarieties $X[n,k](U^+)$, where
$U$ ranges over all subsets of $N$ of cardinality $n-k$.  Fulton and
MacPherson prove that these subvarieties are all disjoint (\cite{FM}).
For $S' = U^+$ where $U \subset N$, $|U| = n-k$, we define
$X[n,k+1](U^+)$ to be the exceptional divisor in $X[n,k+1]$ lying over
$X[n,k](U^+)$. For $S' \in N^+$ not of this form, we define
$X[n,k+1](S')$ to be the proper transform of $X[n,k](S')$.

For each $U \subset N$ of cardinality $|U|=n-k$ we
have the following pullback diagram:

\[\xymatrix{
X[n,k+1](U^+) \ar@{=}[r] & P(\mathcal{N}_{U}) \ar[r] \ar[d] & X[n,k+1]
\ar[d] \\
& X[n,k](U^+) \ar[r] & X[n,k],}\]
where $\mathcal{N}_U=N_{X[n,k](U^+)}X[n,k]$ is the normal bundle of
$X[n,k](U^+)$ in $X[n,k]$ and $X[n,k+1](U^+)=P(\mathcal{N}_U)$
is the exceptional divisor of the blowup. We write $l_U =
c_1(\mathcal{O}_{\mathcal{N}_U}(1))$.

To complete the construction of $X[n+1] = X[n,n]$, we define for $S' \in
N^+$, $X[n+1](S)
= X[n,n](S')$ and $\pi : X[n+1] \ra X^n$ to be the composition
\[\xymatrix{
X[n+1] = X[n,n] \ar[rr]^{\rho_1 \circ \cdots \circ \rho_n} & & X[n,0] =
X[n] \times X \ar[rr]^{\ \ \ \ \ \ \ \ \pi_n \times id_X} & & X^n
}\]

For convenience of notation, we define $X[n,i](S_1, \ldots, S_k) =
X[n,i](S_1) \cap \cdots \cap X[n,i](S_k)$.

\begin{theorem} \label{fm_blowup}
Let $\emptyset \neq S \subset N \neq 2$, $|S| = i$. Choose $a \in N$. Then:
\begin{enumerate}
\item \label{no_change}
The morphisms $X[n,n-1](\{a\}^+) \ra X[n]$ and $X[n,n-s](S^+) \ra
X[n,0](S^+) \cong D(S)$ are isomorphisms for $|S| \geq 2$.
\item \label{fm_bundle_part}
The morphism $X[n,1](N^+) \ra X[n,0](N^+)$ is a projective bundle
morphism of relative dimension $d$.
\item \label{fm_blowup_part}
The morphism $X[n,i+1](N^+) \ra X[n,i](N^+)$ for $1 \leq i \leq n-1$ is
a blowup along the union of subvarieties $X[n,i](N^+, S^+)$ for $S
\subset N$, $|S| = n - i$. 
\item \label{fm_smaller_no_change} For $S$ as above, the morphism
$X[n,i](N^+, S^+) \ra X[n,0](N^+, S^+)$ is an isomorphism.
\end{enumerate}
\end{theorem}

\begin{proof}
Parts \ref{no_change}, \ref{fm_bundle_part}, and \ref{fm_blowup_part}
follow from \cite{FM} proposition 3.5.  For part
\ref{fm_smaller_no_change}, we have by part \ref{no_change},
$X[n,i](S^+) \ra X[n,0](S^+)$ is an isomorphism. Consequently, when we
restrict to $X[n,i](N^+, S^+)$, we get an isomorphism $X[n,i](N^+, S^+)
\cong X[n,0](N^+, S^+)$.
\end{proof}

\subsection{Functors of points related to the Fulton-Macpherson space}

It will be useful for us to have a description of the functors
represented by the varieties $X[n,i]$ and $X[n,i](S_1, \ldots, S_k)$.
Suppose $H$ is a variety and $h : H \ra X^{n+1}$ is a morphism. For a
subset $S \subset N$, we use the notation $h_S$ to denote the
composition of $h$ with the projection $X^{N^+} \ra X^S$, and we write
$h_a$ for $h_{\{a\}}$. Let $\Delta_S \subset X^S$ denote the (small)
diagonal, and $\cI_S$ the ideal sheaf of $\Delta_S$ in $X^S$. Note that
for $S_1 \subset S_2$, the projection $X^{S_2} \ra X^{S_1}$ induces a
morphism $h_{S_1}^* \cI_{S_1} \ra h_{S_2}^* \cI_{S_2}$.

\begin{defn} \label{comp_def}
Let $H$ and $h$ be as above. A \emph{screen} for $h$ and $S \subset N^+$ is 
an invertible quotient $\phi_S : h^*\cI_S \ra \cL_S$. A collection
of screens $\phi_{S_1}, \ldots, \phi_{S_k}$ is \emph{compatible} 
if whenever $S_i \subset
S_j$, there is a unique morphism $\cL_{S_i} \ra \cL_{S_j}$ which makes the
following diagram commute:
\[\xymatrix{
h_{S_i}^*\cI_{S_i} \ar[rr] \ar[d] & & \cL_{S_i} \ar[d] \\
h_{S_j}^*\cI_{S_j} \ar[rr] & & \cL_{S_j} }\]
\end{defn}

\begin{defn}
A subset $S \subset N^+$ satisfies \emph{property} $P_i$,
($S \in P_i$) if either $S \subset N, |S| \geq 2$ or $S = T^+,
|T| > n-i$.
\end{defn}

\begin{defn}
We define the functor $\mc X[n,i]$ from the category of schemes to the
category of sets by setting $\mc X[n,i](H)$ to be the set of pairs
\[\Big((h : H \ra X^{N^+}), \{\phi_T : h_T^* \cI_T \ra \cL_T\}_{T \in
P_i}\Big),\]
 such that the $\phi_T$'s form a compatible collection of screens. We
define the subfunctor $\mc X[n,i](S_1, \ldots, S_k)$ by setting $\mc
X[n,i](S_1, \ldots S_k)(H)$ to be the subset of $\mc X[n,i]$ such that
whenever $T \in P_i$ with $|T \cap S_j| \geq 2$ and $T \not\subset S_j$
for some $j$, the compatibility morphism $\cL_{T \cap S_j} \ra \cL_T$ is
zero.  
\end{defn}

The following theorem is useful not only for understanding the iterative
construction of the space $T_{d,n}$, but also for the applications to
the Fulton-MacPherson configuration space in Section \ref{applications}.

\begin{theorem}[\cite{FM}] \label{fm_functor}
The functors $\mc X[n,i]$ and $\mc X[n,i](S_1, \ldots, S_k)$ are
represented by the varieties $X[n,i]$ and $X[n,i](S_1, \ldots, S_k)$
respectively.
\end{theorem}

\subsection{Inductive construction of $T_{d,n}$}

We now present an inductive construction of $T_{d,n}$ as a sequence of
blowups of a projective bundle over $T_{d,n-1}$. This allows us to give
an explicit inductive presentation of its Chow groups, its Chow motive
and in the next section, a description of its Poincar\'e polynomial.  

\begin{theorem} \label{main_blowup}
There is a sequence of smooth varieties $F_{d,n}^{i}$, for $0 \leq i \leq n$,
with subvarieties $F_{d,n}^i(T)$ indexed by $T \subsetneq N^+$, $|T| \geq 2$ and
morphisms $b_i: F_{d,n}^{i+1} \ra F_{d,n}^i$ such that:
\begin{enumerate}
\item 
$F_{d,n}^0 = T_{d,n}$ and $F_{d,n}^n = T_{d,n+1}$,
\item
the morphism $b_0 : F_{d,n}^1 \ra F_{d,n}^0 = T_{d,n}$ is a projective
bundle morphism of relative dimension $d$,
\item \label{blowup}
the morphisms $b_i$ for $1 \leq i \leq n-1$ are blowups along the union
of the subvarieties $F_{d,n}^i(S^+)$, $|S| = n-i$,
\item \label{products}
$F_{d,n}^i(S^+) \cong T_{d,s} \times T_{d,n-s+1}$ where $s = |S| = n-i$
if $i \neq n-1$,
\item \label{final_noproduct}
$F_{d,n}^{n-1}(\{a\}^+) \cong T_{d,n}$.
\end{enumerate}
\end{theorem}

This theorem closely parallels the above situation for the Fulton
MacPherson configuration space $X[n+1]$ from $X[n] \times X$, and in
fact we  derive it mainly as a consequence of a careful analysis
of certain aspects of this construction. Although very similar in
structure to the construction of Keel \cite{Keel} in the case $d =
1$, we note that our construction is different. For example, the analog of
our (nontrivial) projective bundle $b_0$ in Keel's construction is 
always a trivial vector  bundle.

The varieties $F_{d,n}^i$ and $F_{d,n}(S)$ are defined as follows:

\begin{defn} \label{fdn_def}
Let $X$ be a smooth projective $d$-dimensional variety, and choose $x
\in X(k)$. We abuse notation slightly (by using the symbol $\pi$ in two
different ways) and let $\pi : X[n,i](N^+) \ra X$ be the composition
\[\xymatrix{ X[n,i](N^+) \ar[r] & X[n,i] \ar[rr]^{\rho_1 \cdots \rho_i}
& & X[n] \times X \ar[rr]^{\pi \times id_X} & & X^{n+1} \ar[r] & X }\]
where the final arrow is any of the projections (all give the same
result). We define $F_{d,n}^i$ to be $\pi^{-1}(x) = X[n,i](N^+) \times_X
x$ and $F_{d,n}^i(S) = X[n,i](N^+, S^+) \times_X x$. We also define
$T_{d,n}(S) = F_{d,n}^0(S^+) \subset F_{d,n}^0 = T_{d,n}$.
\end{defn}

The proof of Theorem \ref{main_blowup} follows from the
following proposition combined with Theorem \ref{fm_blowup}.

\begin{proposition} \label{triv_tang}
Suppose $X$ is a smooth variety with trivial tangent bundle, and let
$S \subset N$, $|S| = s$. Then if $0 \leq i \leq s$, we have an
isomorphism $X[n,n-s+i](S^+) = X[n-s+1] \times F_{d,s}^i$
commuting with the natural projections to $X[n-s+1]$.
\end{proposition}

In particular, in proving this proposition, it immediately follows that
Definitions \ref{fdn_def} and \ref{def_tdn} do not depend on the variety
$X$ or on the point $x$, since the tangent bundle is locally trivial.
Note that it follows that for $X$ with trivial tangent bundle, we have
$X[n,n-s+i](S^+, T^+) \cong X[n-s+1] \times F_{d,n}^i(T^+)$ for $T
\subsetneq S$ by restricting the above isomorphism.

\begin{remark} \label{chow_classes}
We see by this reasoning that in definition \ref{fdn_def}, the maps $\pi
: X[n,i](N^+) \ra X$ and $\pi|_{X[n,i](N^+, S^+)}$ locally have a
product structure as above, and so in particular, $[F_{d,n}(S)] =
\iota_x ^! [X[n,i](N^+, S^+)]$ in the Chow groups of $X[n,i](N^+)$ and
$F_{d,n}$ respectively, where $\iota_x$ is the inclusion of the point
$x$ into $X$. Since the class $\big[X[n,i](N^+, S^+)\big]$ on
$X[n,i](N^+)$ can also be seen as the Gysin pullback
$j^!\big[X[n,i](S^+)\big]$ where $j : X[n,i](N^+) \hra X[n,i]$ is the
natural embedding, we may write $[F_{d,n}(S^+)] = \iota^!j^!
\big[X[n,i](S^+)\big]$, and it makes sense therefore to abuse notation
and formally define $[F_{d,n}(N^+)] = \iota^!j^! \big[X[n,i](N^+)\big]$.
We also set $[T_{d,n}(S)] = [F_{d,n}^0(S^+)]$.
\end{remark}

\begin{proof}[proof of proposition \ref{triv_tang}]
By Theorem \ref{fm_functor}, we may translate a morphism to
$X[n,n-s+i](S^+)$ to a collection of screen data. Suppose $f : H
\ra X[n,n-s+i](S^+)$ corresponds to a collection of screens \[\Big((h :
H \ra X^{N^+}), \{\phi_T : h_T^* \cI_T \ra \cL_T\}_{T \in P_i}\Big).\] 
Note that
if $T \subset S$ or $T = S^+$, the function $h_T$ maps entirely into
the diagonal $\Delta_T \subset X^T$, and so by factoring $h_T$ through
$\Delta_T \cong X$, we may identify $h_T^* \cI_T$ with the pullback of
the conormal sheaf $h_{\{a\}} \cI_T/\cI_T^2 = \big((TX)^T/TX\big)^*$,
where $a \in T$. Note that we have $h_{\{a\}} = h_{\{b\}}$ for any $b
\in S^+$, and so it suffices to choose $a \in S^+$. The screen for $T$
therefore translates to specifying a $1$-dimensional sub-vector bundle
$i_T : L_T \hra (TX)^T/TX$ compatible with respect to the natural
projections $(TX)^T/TX \ra (TX)^{T'}/TX$ for $T' \subset T$. Since $V$
is trivial, we may write $TX = V_X$ where $V$ is a vector space over
$k$ of dimension $d$. Let $\mc F_{d,s}(H)$ denote the set of collections
of compatible line sub-vector bundles of the form $\{i_T : L_T \hra
V^T/V\}_{T \subset S \text{ or } T = S^+}.$ It is not hard to check
from this description that $\mc F_{d,s}$ in fact defines a functor which
is represented by $F_{d,s}$.

To complete the proof, we note that by the above, a morphism to
$X[n,n-s+i](S^+)$ gives a morphism to $F_{d,s}^i$. We may also obtain a
morphism to $X[n-s+1]$ by choosing a point $a \in S^+$ and keeping only
the screens for subsets $T \subset (N \setminus S) \cup \{a\}$.  This
gives a morphism $X[n,n-s+i](S^+) \ra X[n-s+1] \times F_{d,s}$. We may
obtain an inverse morphism by describing how to take the screens for $T
\subset (N \setminus S) \cup \{a\}$ together with the screens for $T
\subset S^+$ and define screens for the remaining $T \in P_{n-s+1}$.
This is done in a way similar to the proof of Theorem \ref{fm_blowup},
part \ref{no_change} and we leave the detailed proof to the reader. We
mention as a guide that for $T$ of the form $T = U \cup R$ where $U
\subset S^+$ and $R \cap S^+ = \emptyset$, we set $\cL_T = \cL_{R \cup
\{a\}}$ if $R \neq \emptyset$, and otherwise $\cL_T = L_T^*$ if $T
\subset S^+$.
\end{proof}

\begin{proof}[proof of Theorem \ref{main_blowup}]
Note that we may obtain the
product decomposition $T_{d,n}(S) \cong T_{d,s} \times
T_{d,n-s+1}$, where $s = |S|$ by letting $X$ be a variety with trivial
tangent bundle and considering the commutative diagram where the upper
square is a pullback:

\[\xymatrix{ 
X[n,n-s+i](S^+, N^+) \ar[rr] \ar[d]_{\text{fibers $\cong F_{d,s}^i$}} & &
X[n,n-s+i](S^+) \ar[d]^{\text{fibers $\cong F_{d,s}^i$}} \\
X[n-s+i,0]\Big(\big((N \setminus S) \cup \{a\}\big)^+\Big)
\ar[d]_{\text{fibers $\cong T_{d,n-s+1}$}} \ar[rr] & & X[n-s+1] \ar[d] \\
X \ar@{=}[rr] & & X }\]
\end{proof}

\subsection{Basic properties of $T_{d,n}$}

We now give an example of these spaces for ``minimal'' values of $n$:

\begin{proposition} \label{base_tdn}
We have isomorphisms $T_{1,3} \cong \bP^1$ and for $d > 1$, $T_{d,2} \cong
\bP^{d-1}$. Under these identifications, $[T_{1,3}(\{1,2,3\})] =
\cO_{\bP^1}(-1)$, and $[T_{d,2}(\{1,2\})] = \cO_{\bP^{d-1}}(-1)$.
\end{proposition}

\begin{proof}
First consider the case of $T_{d,2}$. We know that $\A^d[2] =
Bl_{\Delta}(\A^d \times \A^d)$, and that $D(\{1,2\})$ is the
exceptional divisor of the blowup \cite{FM}. Therefore
$$D(\{1,2\}) = P(N_{\Delta} (\A^d \times \A^d)) \cong \bP^{d-1} \times
\A^d.$$ In particular, $T_{d,2} \cong \bP^{d-1}$ as claimed.

For the case of $T_{1,3}$, we note that $\A^1[3] = Bl_{\Delta}
(\A^1)^3$, where $\Delta$ is the small diagonal with exceptional
divisor $D(\{1, 2, 3\})$.
\end{proof}

\begin{corollary}\label{blowupcor1}
$T_{d,n}$ is a smooth projective rational variety of dimension
$dn-d-1$.
\end{corollary}
\begin{proof}
This follows from induction on $n$, with base case proven in 
Proposition \ref{base_tdn} and inductive step given by
 Theorem \ref{main_blowup}.
\end{proof}

\begin{proposition} \label{t1n}
$T_{1,n} \cong \ov{M}_{0,n+1}$.
\end{proposition}
\begin{proof}
Since each $n$ pointed stable $1$-RTPS is exactly an $(n+1)$-pointed
stable curve (where the root hyperplane is identified with the
$(n+1)$st marking), the family $T_{d,n}^+ \ra T_{d,n}$ gives a morphism
$T_{d,n} \ra \ov{M}_{0,n}$ which is bijective on $\ov{k}$ points. To show this
is an isomorphism, it suffices to construct an inverse morphism.

To do this, consider the Fulton-MacPherson configuration space $\bP^1[n]$, and
the divisor $D(N)$ on $\bP^1[n]$. By Definition \ref{def_tdn}, we may
identify $T_{1,n}$ with a fiber of the natural morphism $D(N) \ra \bP^1$. By
\cite{MM} (pages 4,5), there is a natural isomorphism $\ov{M}_{0,
n}(\bP^1, 1) \cong \bP^1[n]$. Note that by \cite{FP} (Theorem 2, part
3), $\ov{M}_{0,n}(\bP^1, 1)$ is actually a fine moduli space for stable maps of
degree $1$ to $\bP^1$. The isomorphism from \cite{MM} identifies the
divisor $D(N)$ with the stable maps $f : (C, p_1, \ldots, p_n) \ra \bP^1$ which
take all the marked points to a given point $x \in \bP^1$. The natural morphism
$D(N) \ra \bP^1$ is simply the morphism taking this stable map to $x$. For such
a stable map, the semistable curve $C$ must have the form:

\medskip

\centerline{
\epsfxsize=2.6in
\epsfysize=1.0in
\epsfbox{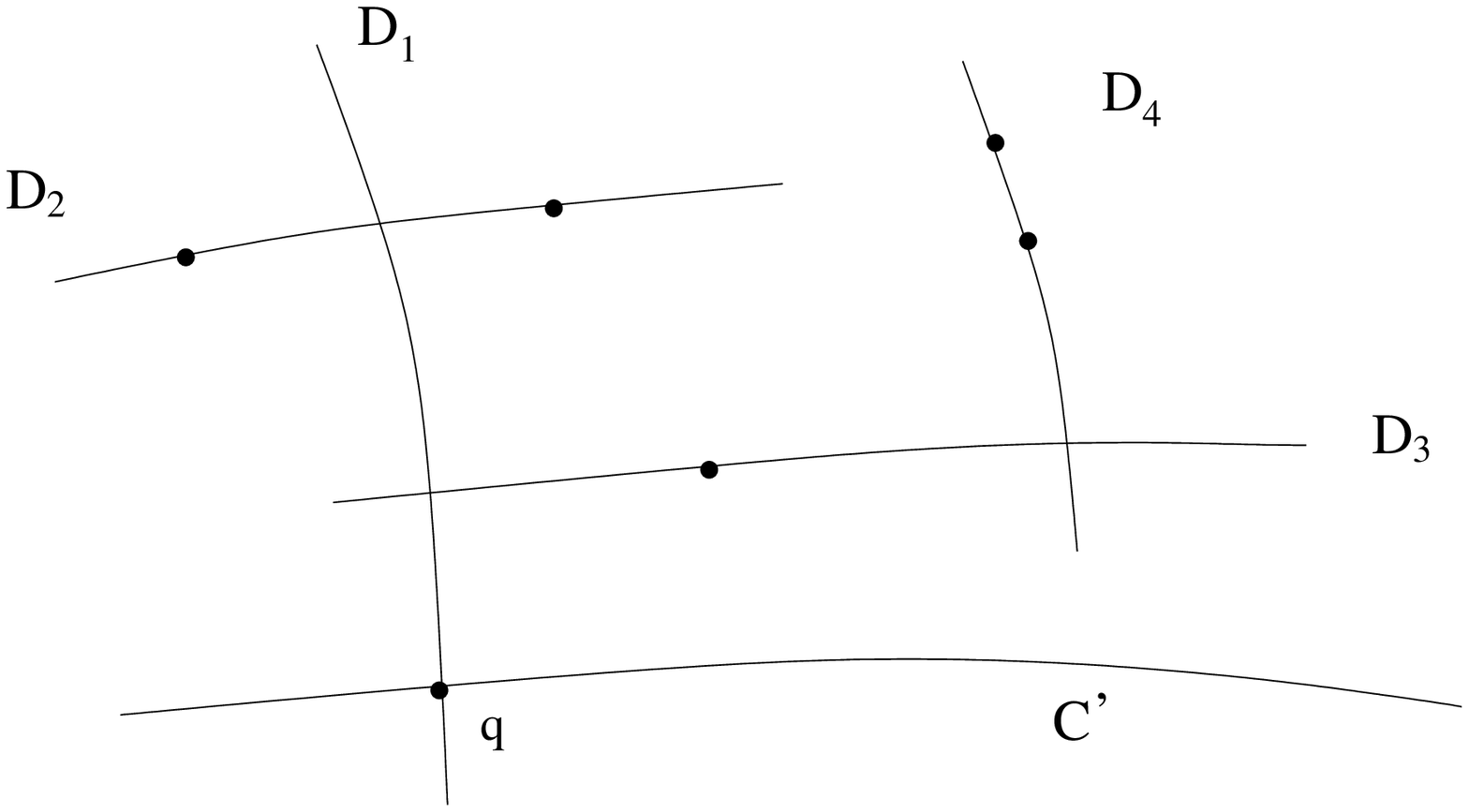}
}

\medskip

Since $\ov{M}_{0,n}(\bP^1, 1)$ is a fine moduli space, one may check
that the fiber over a given point $x \in \bP^1$ is also a fine moduli
space for $(n+1)$-pointed stable curves.  To see this, suppose we have a
stable map $f : (C, p_1, \ldots, p_n) \ra \bP^1$ in this fiber, where
the semistable curve $C$ has irreducible components $C', D_1, \ldots,
D_r$, with $f_*([C']) = 1$, $f_*([D_i]) = 0$. Since none of the marked
points $p_i$ may lie on $C'$, the  curve remaining after forgetting
$C'$, composed of the union of the $D_i$, is a $(n+1)$-pointed stable
rational curve where the $(n+1)$st marking is obtained from the point
where $\cup D_i$ intersects $C'$. This gives an isomorphism of the fiber
over $x$ with $\ov{M}_{0, n+1}$ which is inverse to the morphism above.
\end{proof}

\begin{theorem} \label{tree_thm}
The points of $T_{d,n}$ are in one to one correspondence with
isomorphism classes of $n$-pointed stable rooted trees of
$d$-dimensional projective spaces.
\end{theorem}
\begin{proof}
Let $X[n]^+ \ra X[n]$ be the ``universal family'' as described in
\cite{FM}. By base change over $T_{d,n} \ra X[n]$ (included by
choosing a point $x \in X$), we obtain a flat family $T_{d,n}^+ \ra
T_{d,n}$. It follows from the description in \cite{FM} that the fibers
are all $n$ pointed stable $d$-RTPS's, exactly one in each isomorphism
class.
\end{proof}

\subsection{Ample divisor classes on $T_{d,n}$}

Although by \ref{blowupcor1}, we know abstractly that $T_{d,n}$ is a
projective variety, it is often useful to have an explicit
presentation of an ample divisor. We exhibit such a divisor below:

\begin{theorem} \label{ample}
Let $\delta_S = [T_{d,n}(S)]$ in the Chow group of $T_{d,n}$ as in
remark \ref{chow_classes}. Then for $S \subset N, |S| \geq 2$, the
divisor classes
$$ \eta_S = \underset{N \supseteq T \supseteq S}{\sum} \delta_T$$
is nef and base point free. Furthermore, any expression of the form
$$ A = \underset{S \subset N, |S| \geq 2}{\sum} c_S \eta_S$$ is
very ample provided $c_S >0$ for all $S$.
\end{theorem}
\begin{proof}
By \cite{FM}, we have for any smooth $d$-dimensional variety $X$, an
embedding $$i : X[n] \hra \prod_{S \subset N, |S| \geq 2}
Bl_{\Delta_S}(X^S),$$ where $\Delta_S \subset X^S$ is the (small)
diagonal. Let $i_S$ be the morphism $i$ composed with the projection map
onto the factor $Bl_{\Delta_S}(X^S)$. If $A_S$ is an very ample class on
$Bl_{\Delta_S}(X^S)$, then it follows immediately that $i_S^*(A_S)$ is
nef and base point free (since it is the pullback of a very ample
divisor), and that $\sum c_S i_S^*(A_S)$ is very ample on $X[n]$ if each
$c_S > 0$ (since it is the pullback of a very ample divisor via an
embedding).

Let $B_S'$ be a very ample divisor class on $X^S$, and let $I_S$ be
the ideal sheaf of $\Delta_S$ in $\cO_{X^S}$. Let $f_S :
Bl_{\Delta_S}(X^S) \ra X^S$ be the natural projection. Then
$f_S^{-1}(I_S)$ is an invertible sheaf and the divisor class $\lambda
f_S^*B_S + c_1(f^{-1}(I_S))$ is very ample for some $\lambda > 0$ (see
\cite{Har}, Proposition 7.10(b) and the proof of 7.13(a)). Let $B_S =
\lambda B_S'$. Then
$$\alpha_S = i_S^*\Big(f_S^*B_S + c_1\big(f^{-1}(I_S)\big)\Big)$$ is
nef and base point free. Fix for the remainder of the proof integers
$c_S > 0$. Then
$$\alpha = \underset{S \subset N, |S| \geq 2}{\sum}
i_S^*\Big(f_S^*B_S + c_1\big(f^{-1}(I_S)\big)\Big)$$ is very ample on
$X[n]$.

Now consider the embedding $T_{d,n} \overset{j}{\hra} X[n]$, and let
$j_S$ be the composition with the projection to
$Bl_{\Delta_S}(X^S)$. Let $\eta_S = j^*\alpha_S$, and $A =
j^*\alpha$. Since $j$ is an embedding, $\eta_S$ is nef and base point
free and $A$ is very ample. We will be done once we show that the
divisors $\eta_S$ have the desired form. To begin, we may rewrite
$\eta_S$ in the following way:
\begin{align*}
\eta_S &= j^*\alpha_S = f_S^*B_S + c_1\big(f_S^{-1}(I_S)\big) \\
&= j^*i_S^*f_S^*B_S + j^* c_1\big(i_S^* f_S^{-1}(I_S)\big) \\
&= (f_S \circ i_S \circ j)^* B_S - j^* c_1\big(i_S^*
f_S^{-1}(I_S)\big) \\
\end{align*}

Examining the second term, we see by lemma \ref{ht1}, that $i_S^*
f_S^{-1}(I_S) = (f_S \circ i_S)^{-1}(I_S)$, and by \cite{FM}, page
203, if we let $I(D(S))$ be the ideal sheaf of $D(S) \subset X[n]$, we
have $(f_S \circ i_S)^{-1}(I_S) = \prod_{T \supset S}
I(D(T))$. Consequently, taking first Chern classes and applying $j^*$,
we have:

$$j^* c_1\big(i_S^* f_S^{-1}(I_S)\big) = j^* \sum_{T \supset S}
c_1(I(D(T))) = j^* \sum_{T \supset S} -[D(T)] = - \sum_{T \supset S}
\delta_T.$$

On the other hand, looking at the first term, $(f_S \circ i_S \circ
j)^* B_S$, we see that since the morphism $(f_S \circ i_S \circ j)$
factors through a morphism to $Spec(k)$, whose Picard group is the
zero group, this pullback must in fact vanish. Therefore we have
$$\eta_S = \underset{S \subset T \subset N}{\sum} \delta_T,$$
as desired.
\end{proof}

\subsection{A relative version of $T_{d,n}$}

It will be useful in what follows to have a relative version of
construction of $T_{d,n}$. This follows without  much technical
difficulty and we leave some of the routine verifications to the reader.

\begin{defn} \label{sheaf_functor}
Let $V$ be a rank $d$ vector bundle over a scheme $X$. We define the
functor $\mc T_{V,n}$ from the category $(Sch/X)^{op}$ of $X$-schemes to
the category of sets as follows. For $h : H \ra X$, we define
$\mc T_{X,n}(H)$ to be the set of collections $\{\phi_T : h^*V \ra
\cL_T\}_{T \in N, |T| \geq 2}$, such that the $\phi_T$'s form a
compatible collection of screens (in the same sense as in
\ref{comp_def}). Note that there is a canonical morphism (natural
transformation) from $\mc T_{V,n}$ to $X$.
\end{defn}

By setting $T_{V,1} = X$, $T_{V, 2} = Proj_X(Sym^\bullet V)$, we may
inductively define varieties $T_{V,n}, F_{V,n}^i, F_{V,n}(S^+)$ such
that the following theorem holds:

\begin{theorem} \label{rel_blowup}
There is a sequence of schemes $F_{V,n}^{i}$, for $0 \leq i \leq n$,
with subschemes $F_{V,n}^i(T)$ indexed by $T \subsetneq N^+$, $|T| \geq
2$ and morphisms $b_i: F_{V,n}^{i+1} \ra F_{V,n}^i$ such that:
\begin{enumerate}
\item \label{bundle_part}
$F_{V,n}^0 = T_{V,n}$ and $F_{V,n}^n = T_{V,n+1}$,
\item
the morphism $b_0 : F_{V,n}^1 \ra F_{V,n}^0 = T_{V,n}$ is a projective
bundle morphism of relative dimension $d$,
\item
the morphisms $b_i$ for $1 \leq i \leq n-1$ are blowups along the union
of the subschemes $F_{V,n}^i(S^+)$, and $F_{V,n}^i(S^+) \cong T_{V,s}
\times_X T_{V,n-s+1}$ if $i \neq n-1$ and $F_{V,n}^{n-1}(\{a\}^+) \cong
T_{V,n}$,
\item
$\mc T_{V,n}$ is represented by $T_{V,n}$.
\end{enumerate}
\end{theorem}

The proof of this theorem takes the following steps. First, we see
that it holds in the case that $V$ is a trivial bundle by noting that
the inductive construction of the space $T_{d,n+1}$ from $T_{d,n}$ by
taking a projective bundle and blowing up may be fibered with a scheme
$X$ to give an inductive construction of $T_{d,n+1} \times X =
T_{V,n+1}$ from $T_{d,n} \times X = T_{V,n}$. For a general bundle, we
may define the functors $\mc F_{V,n}^i(S^+)$ by emulating the definition
of $\mc X[n,i](N^+, S^+)$, and note that they are locally represented
subschemes over subsets where $V$ is trivial, and that these subschemes
glue to give a closed subscheme $F_{V,n}^i(S^+)$. We define
$F_{V,n}^{i+1}$ to be the blowup along the subschemes $F_{V,n}^i(S^+)$
where $|S| = i$. Since these have the correct functorial description
locally, we may glue and conclude that $F_{V,n}^{i+1}$ represents the
functor $\mc F_{V,n}^{i+1}$.

We note the following lemmas which will be useful in Section
\ref{applications}:

\begin{lemma} \label{prod_tdn}
Let $V$ be a vector bundle on $X$, and let $\pi : T_{V,n} \ra X$ be the
natural projection. Then $T_{V,n} \times_X T_{V,m} = T_{\pi^* T_{V,n},
m}$.
\end{lemma}

\begin{lemma} \label{triv_tdn}
Let $X$ be any scheme, and let $V$ be a trivial vector bundle of rank
$d$ on $X$. 
Then $T_{d,n} \times X \cong T_{V,n}$. In particular, if $X = Spec(k)$
then $T_{V,n} \cong T_{d,n}$.
\end{lemma}

\begin{lemma} \label{subset_tdn}
Let $X$ be a smooth variety and let $D(S) \subset X[n]$ be the divisor
on the Fulton-MacPherson configuration space described in the beginning
of the section. If $\pi : D(S) \ra X[n-|S|+1] = X[(N \setminus S) \cup
\{a\}]$ is the natural morphism from projecting with respect to the
subset $(N \setminus S) \cup \{a\}$ where $a \in S$, and $f : X[(N
\setminus S) \cup \{a\}] \ra X$ the projection with respect to $a$, then
there is a natural isomorphism $D(S) = T_{f^* \Omega_X, s}$.
\end{lemma}

The proofs of these are elementary and follow from an examination of
the functorial descriptions of the spaces involved.

\begin{remark} \label{point_drop} 
The morphism $T_{d,n+1} \ra T_{d,n}$ obtained by composing the morphisms
$b_i$ of theorem \ref{rel_blowup} is given functorially by dropping all
screens for subsets $S \subset N^+$ which contain the $(n+1)$st marking.
In the future, we denote this morphism by $\pi_{n+1}$. We may similarly
define morphism $\pi_i$ for any $i \in N^+$ by dropping the $i$th
marking.
\end{remark}

\section{Inductive presentations of Chow groups and motives}
\label{applications}

We consider the Chow group of a variety $X$ as a graded abelian group
$A(X) = A_*(X)$. We use the following general conventions. For a graded
abelian group $M = \oplus_{i \in \Z} M_i$, we set $M(n)$ to be the group
with grading shifted so that $(M(n))_i = M_{n-i}$. The grading on the
tensor product $M \otimes N = M \otimes_\Z N$ is given by $(M \otimes
N)_n = \oplus_{i+j = n} M_i \otimes N_j$. We write $\Z$ for the graded
abelian group with the integers in degree $0$ and zero in all other
degrees. 

\begin{remark}
All the constructions used in this section are motivic: in
other words, if the reader prefers, they may interpret $A(X)$ as $M(X)$,
the Chow motive of $X$ (as in \cite{Man:Mot}), $M(i)$ to be $M$ twisted
$i$ times with the Lefschetz motive, and $\Z$ to be the motive of
$Spec(k)$.
\end{remark}

In this notation we have the following well known facts:

\begin{lemma} \label{proj_bundle_chow}
For $V$ a vector bundle of rank $d$ on $X$, 
\[A(P(V)) = \bigoplus_{i = 0}^{d-1}A(X)(i) = A(X) \otimes (\bigoplus_{i
= 0}^{d-1} \Z(i))= A(X) \otimes A(\bP^{d-1}).\]
\end{lemma}
\begin{proof}
\cite{Ful:IT} for Chow groups, \cite{Man:Mot} for Chow motives.
\end{proof}

\begin{lemma} \label{blowup_chow}
For $Z \hra X$ a regularly embedded subvariety of codimension $d$,
$$A(Bl_Z X) = A(X) \oplus \bigoplus_{i = 1}^{d-1}A(Z)(i) = A(X) \oplus \Big(A(Z)
\otimes A(\bP^{d-2}(1))\Big).$$
\end{lemma}
\begin{proof}
\cite{Man:Mot}
\end{proof}

One technical difficulty which makes the computation of Chow groups more
difficult than the computation of cohomology is the fact that the Chow
group of a product $X \times Y$ is not easily expressible in terms of
the Chow groups of $X$ and $Y$. Philosophically, we show in this section
that products (and certain fiber bundles) with $T_{d,n}$ as one of the
factors are not subject to this difficulty.

\subsection{The Chow groups and motives of $T_{d,n}$} \label{tdn_chow}

\begin{theorem} \label{tvn_chow_groups}
Let $V/X$ be a vector bundle of rank $d$. Then
\begin{multline*} 
A(T_{V,n+1}) =
\left(
\bigoplus_{j=0}^{d} 
A(T_{V,n})(j)
\right) 
\oplus \\
\left(
\bigoplus_{j=1}^d 
\bigoplus_{S \subsetneq N, |S| \geq 2} 
A(T_{V,|S|} \times T_{V,n-|S|+1})(j)
\right)
\oplus
\left(
\bigoplus_{j=1}^{d-1}
\bigoplus_{a \in N}
A(T_{V,n})(j)
\right)
\end{multline*}
\end{theorem}
\begin{proof}
This is an immediate consequence of Theorem \ref{rel_blowup} together
with Lemmas \ref{proj_bundle_chow} and \ref{blowup_chow}.
\end{proof}

Using Lemma \ref{triv_tdn} to identify $T_{d,n} \times B = T_{\cO_B^d,
n}$ for a variety $B$, we obtain the following corollary:

\begin{corollary} \label{tdn_chow_groups}
Let $V/X$ be a vector bundle of rank $d$. Then
\begin{multline*} 
A(T_{d,n+1} \times B) =
\left(
\bigoplus_{j=0}^{d} 
A(T_{d,n} \times B)(j)
\right) 
\oplus \\
\left(
\bigoplus_{j=1}^d 
\bigoplus_{S \subsetneq N, |S| \geq 2} 
A(T_{d,|S|} \times T_{d,n-|S|+1} \times B)(j)
\right)
\oplus
\left(
\bigoplus_{j=1}^{d-1}
\bigoplus_{a \in N}
A(T_{d,n} \times B)(j)
\right)
\end{multline*}
\end{corollary}

%

In particular, setting $B = Spec(k)$, we obtain a presentation for the
Chow groups of $T_{d,n+1}$ in terms of the Chow rings of varieties of
the form $T_{d,t} \times B$ for various values of $t<n$. These terms
may be successively reduced using corollary \ref{tdn_chow_groups},
eventually using the fact that $A(T_{1,3}) = A(\bP^1) = \Z \oplus
\Z(1)$ or the fact $A(T_{d,2}) \cong A(\bP^{d-1}) = \oplus_{i =
0}^{d-1} \Z(i)$ from Lemma \ref{base_tdn}.

\subsection{The Chow groups and motives of the Fulton-MacPherson configuration space}

\noindent
Fulton and MacPherson have given a compactification $X[n]$ of the
moduli space of $n$ distinct points on a smooth variety $X$ as well as
a presentation of the intersection ring of the space $X[n]$.  They
pose the following problem ([FM/Ann94], page $189$):

\smallskip

\begin{quote} \label{chowFM}
"It would be interesting to find an explicit basis for the Chow groups
of $\mathbf{P}^m[n]$, preferably simple with respect to the intersection
pairings, as Keel has done in the case $m=1$."
\end{quote}

\smallskip

In this section, we present an inductive presentation of the Chow groups
and motives of the spaces $X[n]$ which parallels Keel's presentation in
\cite{Keel}.

Let us begin by making an elementary observation. The blowup
construction of the Fulton-MacPherson space described in Section
\ref{construction} gives $X[n+1]$ as a composition of blowups of
$X[n]\times X$. In the same way, we easily obtain $X[n+1] \times B$ as a
composition of blowups of $X[n] \times X\times B$ for an arbitrary
variety $B$.  Together with the identifications from Theorem
\ref{fm_blowup}, part \ref{no_change}, this yields the following
inductive presentation:

\begin{theorem} \label{fm_chow_groups}
Let $X$ be a smooth variety. Then
\begin{multline*} 
A(X[n+1] \times B) = 
A(X[n] \times X\times B)  \oplus \\
\left(
\bigoplus_{j=1}^d 
\bigoplus_{S \subsetneq N, |S| \geq 2} 
A(D(S) \times B)(j)
\right)
\oplus
\left(
\bigoplus_{j=1}^{d-1}
\bigoplus_{a \in N}
A(X[n] \times B)(j)
\right)
\end{multline*}
\end{theorem}

As before, the symbol $A$ can stand either for the Chow group or the
Chow motive.  In particular, setting $B = Spec(k)$, we obtain a
presentation for the Chow groups (or motives) of $X[n+1]$ in terms of
the Chow groups (or motives) of the varieties $D(S)$ and $X[n] \times
X$. From Lemma \ref{subset_tdn}, we have an isomorphism $D(S) = T_{V,
s}$, for a vector bundle $V$ of rank $d$ on $X[n-s+1]$. These terms may
be successively reduced using theorems \ref{tvn_chow_groups} and
\ref{fm_chow_groups}, eventually yielding an answer in terms only of the
Chow groups of $X^T$ for $T = 1, \ldots, n$. In general there is no
known formula for the Chow groups of $X^T$ in terms of the Chow groups
of $X$, however, if $X$ has a cellular decomposition (see definition
\ref{cell_def}), then it is true that $A(X^T) = \otimes_{a \in T} A(X)$.

\section{Betti numbers and Poincar\'e polynomials} \label{poincare}

In this section, we analyze generating functions for the Betti numbers
and Poincar\'{e} polynomials of the varieties $T_{d,n}$ in the case when
$k = \C$. By Corollary \ref{tdn_HI}, the Betti numbers coincide with the
ranks of the Chow groups, and therefore, since the presentation of the
Chow groups of $T_{d,n}$ is independent of the underlying field $k$,
these determine the Chow groups in general. A recursive description of
these polynomials was given in \cite{FM} for the spaces $X[n]$. In
\cite{Man:GF}, Manin relates the Poincar\'e polynomials of these spaces
as well as the polynomials for $\M_{0,n}$ to solutions to certain
differential and functional equations. We apply Manin's ideas here to
obtain similar results for $T_{d,n}$, which specialize to Manin's
original result for $\overline{M}_{0,n+1}$ in the case $d = 1$.

Indeed, the defining equations for the generating functions of $T_{d,n}$
described here are identical to those discussed by Manin (Theorem 0.4.1)
in his analysis of the Poincar\'{e} polynomials of $X[n]$.  In our case,
we recover these equations from the explicit blowup construction of
$T_{d,n}$, just as one can recover Theorem 0.3.1 of Manin from the
blowup construction of $\overline{M}_{0,n}$ of Keel.


For a smooth compact variety $Z$, denote  its Poincar\'{e} polynomial by
$P_Z(q)=\sum_j \dim H^j(Z)q^j$.  In particular, put
\[
\kappa_m=P_{\bP^{m-1}}(q)=\frac{q^{2m}-1}{q^2-1}.
\]
Fix $d$, and for $n\geq 2$, denote by $P_n(q)=P_{T_{d,n}}(q)$ the
Poincar\'{e} polynomial of $T_{d,n}$. From corollary \ref{tdn_HI} and
the inductive presentation of the Chow groups of $T_{d,n}$ in section
\ref{tdn_chow}, we have the following recursion for the Poincar\'{e}
polynomials $P_n(q)$.
\begin{equation}
\label{recursionP}
P_{n+1}(q)= (\kappa_{d+1}+nq^2\kappa_{d-1})P_n(q) +
q^2\kappa_d\sum_{{i+j=n+1}\atop{2\leq i\leq n-1}}
{\textstyle{\binom{n}{i}}}P_i(q)P_j(q)
\end{equation}
Defining $P_1(q)=1$, and defining $p_n=p_n(q)=\displaystyle\frac{P_n(q)}{n!}$,
this is equivalent to either of the recursions:
\begin{multline*}
(n+1)p_{n+1}= (\kappa_{d+1}+nq^2\kappa_{d-1})p_n +
q^2\kappa_d\sum_{{i+j=n+1}\atop{2\leq i\leq n-1}} jp_ip_j\\
(n+1)p_{n+1}= (\kappa_{d+1}+nq^2\kappa_{d-1})p_n +
q^2\kappa_d\sum_{{i+j=n+1}\atop{i\geq 1}} jp_ip_j
-q^2\kappa_dnp_n-q^2\kappa_dp_n.
\end{multline*}

The fact that $\kappa_{d+1}=1+q^2\kappa_d$, $q^2\kappa_{d-1}=\kappa_d-1$,
and that $q^2\kappa_d=q^{2d}-1+\kappa_d$ shows that the recursion in
(\ref{recursionP}) can be rewritten as:
\begin{equation}
\label{recursionp}
(n+1)p_{n+1}(q)= (1-nq^{2d})p_n(q) +
q^2\kappa_d\sum_{{i+j=n+1}\atop{i\geq 1}} jp_i(q)p_j(q).
\end{equation}

Consider the following generating function, recalling that $p_1(q)=1$:
\[
\psi(q,t)= t+\sum_{n\geq 2} p_n(q)t^n = \sum_{n\geq 1}p_n(q)t^n.
\]

\begin{theorem}
\label{generating}
$\psi(q,t)$ is the unique root in $t+t^2\Q[q][[t]]$ of the following
functional equation in $t$ with parameter $q$:
\begin{equation}
\label{functional}
\kappa_d(1+\psi)^{q^{2d}} =q^{2d+2}\kappa_d\psi-q^{2d}(q^{2d}-1)t+\kappa_d
\end{equation}
or the following differential equation in $t$ with parameter $q$:
\begin{equation}
\label{differential}
(1+q^{2d}t-q^2\kappa_d\psi)\psi_t =1+\psi.
\end{equation}
\end{theorem}

\begin{proof} First, note that we get (\ref{differential}) from
(\ref{functional}) by differentiating in $t$.
Moreover, we can see that the  equations are equivalent
to the recursion (\ref{recursionP}) or (\ref{recursionp}).
In particular, since $\psi_t(q,t) =  \sum_{n\geq 1} np_n(q)t^{n-1}$,
the $t^n$ term for $n\geq 1$
of the left hand side of the differential equation is
\[
(n+1)p_{n+1}+q^{2d}np_n -q^2\kappa_d \sum_{{i+j=n+1}\atop{i\geq 1,j\geq 1}}jp_jp_i
\]
which is equal to $p_n$ by (\ref{recursionp}), and for $n=0$ is $p_1=1$.
This is exactly the statement of the theorem.
\end{proof}

Fix $d\geq 1$ and define the generating function
$\eta(t)=t+\sum_{n\geq 2}\chi(T_{d,n})$.

\begin{corollary}
$\eta(t)$ is the unique root in $t+t^2\Q[q][[t]]$ of any of the following
equations:
\begin{eqnarray*}
d(1+\eta)\log(1+\eta) &=& (d+1)\eta -t \\
(1+t-d\eta)\eta_t &=& 1+\eta.
\end{eqnarray*}
\end{corollary}

\begin{proof}
Differentiating the first equation gives the second.  The result
follows since the Euler characteristic of a smooth compact variety $Z$ can be
defined by $\chi(Z) = P_z(-1)$ so that $\eta(t)=\psi(-1,t)$.
\end{proof}

\section{Chow Ring of $T_{d,n}$ and pairing between divisors and curves} \label{pairing_section}

In this section we examine the structure of the Chow ring of the space
$T_{d,n}$. Let $k$  be an arbitrary field. We recall the definition of
the divisors $T_{d,n}(S)$ described in Section \ref{construction}. Fix
$d\geq 1, n\geq 2$, and let $\delta_S = [T_{d,n}(S)]$ be the
corresponding cycle class in the Chow group $A_*(T_{d,n})$. We obtain an
explicit description of the Chow ring of $T_{d,n}$ by considering the
Fulton-MacPherson configuration space $\A^d[n]$. Let $i : D(N)
\hookrightarrow \A^d[n]$ be the inclusion of the divisor on $\A^d[N]$
where all the points coincide, and consider the morphism $\pi : D(N) \ra
\A^d$. By Definition \ref{def_tdn}, $T_{d,n}=\pi^{-1}(0)$.  Recall that by
Proposition \ref{triv_tang}, $T_{d,n} \times \A^d\cong D(N)$.

This implies that we may regard $D(N)$ as a (trivial) vector bundle over
$T_{d,n}$, and therefore we have a morphism $\pi : D(N) \ra T_{d,n}$, and
a flat pullback inducing an isomorphism $(\pi^*)^{-1} : A^* D(N) \cong A^*
T_{d,n}$. In particular, one may check from the definitions that if we
let $D(S)$ be the divisor defined in \cite{FM} which we used in Section
\ref{construction} then for $S \neq N$, $$\delta_S = (\pi^*)^{-1} i^!
[D(S)].$$ We may similarly define $\delta_N = (\pi^*)^{-1} i^! [D(N)]$.
For any two distinct elements $a$ and $b \in N$, define $\displaystyle
\Sigma_{ab}:= \sum_{S \subset N \setminus \{ab\}} \delta_{\{ab\} \cup
S}.$

\begin{theorem} \label{Chowring}
$\displaystyle A^*(T_{d,n}) \cong \mathbb{Z}[\{\delta_S\}_{S\subset N, 2 \le |S|
\le n}]/I_{d,n},$ where $I_{d,n}$ is the ideal generated by:
\begin{enumerate}
\item \label{rel1} $\delta_S \cdot \delta_T =0$ for all $S,T \subset
N$, such that $2 \le |S|,|T|$ and $\emptyset \ne S\cap T \subsetneq S,
T$;
\item $(\Sigma_{ij})^d=0$, for all $i$,$j \in N$.
\end{enumerate}
\end{theorem}

\begin{proof}
This follows from immediately from \cite{FM}, and the isomorphism
$A^*(T_{d,n}) \cong A^*(D(N))$ above.
\end{proof}

In the remainder of this section, we state and prove the following 
pairing between  $1$-cycles and the boundary divisors.

\begin{theorem}\label{pairing}
For $T \subsetneq N$, $|T| \ge 2$, define $1$-cycles
$C_T\in A_1(T_{d,n})$ by
$$C_T := \delta_T^{d(|T|-1)-1} \cdot \delta_N^{d(n-|T|)-1}.$$
If $S\subsetneq N$, $|S|\geq 2$, then
\[
\delta_S \cdot C_T = 
\begin{cases} 
(-1)^{d(n-1)} & \text{ if } $S=T$; \\
(-1)^{n-2} & \text{ if $d=1,|S|=2,S\subsetneq T$};\\
0 & \text{ otherwise.}
\end{cases}
\]
\end{theorem}

\begin{corollary}\label{intbasis}
For $d > 1$, the $1$-cycles $C_T$, $|T| \ge 2$ form a
$\Z$-basis for $A_1(T_{d,n})$. In the case $d = 1$,  the
$1$-cycles $C_T$, $|T| \ge 3$ form a $\Z$-basis for
$A_1(T_{1,n})=A_1(\M_{0,n+1})$.
\end{corollary}

\begin{proof}
First note that the set $\{\delta_T\}_{T \subset N, |T|
\ge 3 }$ forms a $\Z$-basis for the codimension $1$-cycles on
$T_{1,n}$, and that the set $\{\delta_T\}_{T \subset N, |T| \ge 2}$ 
forms a $\Z$-basis for the codimension $1$-cycles on
$T_{d,n}$. 

Note that the ring described in the theorem does not depend on the choice of
the base field. The statement which we have to prove is independent of the
choice of $k$, and we may assume without loss of generality that $k = \C$. By
Appendix \ref{tdn_HI}, the space $T_{d,n}$ is an HI space. Since it is a
compact smooth manifold with torsion free cohomology groups, we obtain
Poincar\'e duality induced by the intersection pairing of divisors and
curves. Therefore the $C_T$'s form a dual integer basis to the $D_T$'s, up to
sign.
\end{proof}

In order to prove Theorem \ref{pairing}, we first establish several
identities. For convenience, we denote by
$\delta_i$ the divisor class $\delta_{\{1, 2, \ldots, i\}}$.
The following is an immediate consequence of Theorem \ref{Chowring}.
\begin{lemma}\label{Chowlemma}
For $S,T\subset N$, $i\in S, T$, and $l\in T\setminus S$, we have
$\delta_T\cdot\delta_S=0$ unless $S\subset T$.  In particular, if $l\in
T$, then $\delta_T\cdot\delta_{l-1}=0$ unless $\{1,\ldots, l\}\subset
T$.
\end{lemma}

\begin{lemma} \label{bigvanishing}
For $S \subsetneq N$, $|S|\ge 2$,
$\delta_S\cdot \delta_N^{d(n-|S|)}=0.$  Consequently, if $S\neq N$, $|S|>j$,
then $\delta_S\cdot\delta_N^{d(n-j)-1}=0$.
\end{lemma}
\begin{proof}
We proceed by induction on $n-|S|$.  The base case $n-|S|=0$ holds trivially.
Suppose that the result holds for $n-|S|< k$. Let $S\subset N$, $|S|=n-k$.  
Choose $i\in S$ and $j\not\in S$. By Theorem \ref{Chowring}, if 
$j\in T$ and $i\in S\cap T$, then $\delta_T\delta_S=0$ unless
$S\subsetneq T$.  Moreover, for such $T\neq N$,
$\delta_T\delta_N^{d(n-|S|-1)}=0$ by the inductive hypothesis.
Therefore
\[
0=(\Sigma_{i,j\in T}\delta_T)^d\delta_S\delta_N^{d(n-|S|-1)}
=\delta_N^d\delta_S\delta_N^{d(n-|S|-1)}=\delta_S\delta_N^{d(n-|S|)}
\]
since $(\Sigma_{ij})^d=0$ by Theorem \ref{Chowring}. 
\end{proof}

\begin{lemma} \label{prop1.2}  Given $2\leq j\leq n$,
if $1\le k<i\le j$, $T\subset
\{1,\ldots i\}$, $|T|=i-k$, then
$$\delta_T\cdot \delta_i^{kd+1} \cdot \delta_{i+1}^d \cdots
\delta_j^d \cdot \delta_N^{d(n-j)-1}=0.$$
\end{lemma}

\begin{proof}  Renumbering the elements of $T$, it suffices
to take $\delta_T=\delta_{i-k}$.
We proceed by induction on $k$.  For the base case $k=1$,
we proceed by induction on $j-i$ 
with base case $i=j$. By Lemma \ref{Chowlemma},
if $1,j\in T$, then $\delta_T\delta_{j-1}=0$ unless $\{1,\ldots,j\}\subset T$.
Also note that by Lemma \ref{bigvanishing}, if $T\neq N$ and 
$|T|>j$, then $\delta_T\delta_N^{d(n-j)-1}=0$.  The relation
$(\Sigma_{1j})^d=0$ of
Theorem \ref{Chowring} gives
\[
0=(\Sigma_{1,j\in T} \delta_T)^d \cdot \delta_{j-1} \cdot \delta_j \cdot \delta_N^{d(n-j)-1}
= (\delta_j + \delta_N)^d\cdot \delta_{j-1} \cdot \delta_j \cdot \delta_N^{d(n-j)-1}.
\]
The summands coming from the terms of the
binomial expansion of $(\delta_j + \delta_N)^d$ with positive
degree in $\delta_N$ vanish by Lemma \ref{bigvanishing} for $\delta_S=\delta_j$,
giving the result.

Now suppose that the result holds for integers less than $j-i$, 
If $1,i\in T$, then $\delta_T\delta_{i-1}=0$
unless $\{1,\ldots,i\}\subset T$ by Lemma \ref{Chowlemma}.  Indeed, 
$\delta_T\delta_{i-1}\delta_i\cdots\delta_j=0$ unless 
$T=\delta_i,\ldots,\delta_j$ or $\{1,\ldots, j\}\subsetneq T$. 
The relation $(\Sigma_{1i})^d=0$ of Theorem \ref{Chowring}
gives
\[
0=(\Sigma_{1,i\in T}\delta_T)^d \cdot \delta_{i-1}\cdot \delta_i \cdot
\delta_{i+1}^d \cdots \delta_j^d \cdot \delta_N^{d(n-j)-1}
\]
We see that terms involving $|T|>j$
vanish by  Lemma \ref{bigvanishing}.  Since terms
involving $\delta_T=\delta_{i+1},\ldots,\delta_j$ vanish by the inductive
hypothesis, we have established
the base case $k=1$.

Suppose that the result holds for integers less than $k$, 
and consider the following identity:
\[
0=(\Sigma_{1,i-k+1\in T}\delta_T)^d \cdot \delta_{i-k} \cdot \delta_i^{(k-1)d+1}
\cdot \delta_{i+1}^d \cdots \delta_j^d \cdot \delta_N^{d(n-j)-1}
\]
By Lemma \ref{Chowlemma},
$\delta_T\delta_i\ldots\delta_j=0$ 
unless $\{1,\ldots,i-k+1\}\subset T\subset\{1,\ldots,i\}$, 
$\{1,\ldots,j\}\subset T$, or
 $\delta_T=\delta_i,\ldots,\delta_j$.  
Terms involving 
$\delta_T$, of the first type vanish by the inductive hypothesis, 
and  terms involving
$|T|> j$ are zero  by Lemma \ref{bigvanishing}.
Therefore, nonzero terms can involve only
$\delta_T=\delta_i,\ldots,\delta_j$.  Finally, the terms with contributions
from $\delta_{i+1},\ldots,\delta_j$ vanish by the base case $k=1$.
Hence we are left with
$$
0=\delta_i^d\delta_{i-k} \cdot \delta_i^{(k-1)d+1} \cdot \delta_{i+1}^d \cdots \delta_j^d \cdot \delta_N^{d(n-j)-1},
$$
proving the proposition.                       
\end{proof}

\begin{lemma}\label{squish}
For $j \ge 1$, $\displaystyle \delta_2^d \cdots \delta_j^d \cdot
\delta_N^{d(n-j)-1} =(-1)^{j-1}\delta_N^{d(n-1)-1}.$
\end{lemma}

\begin{proof}
We prove the result by induction on $j$.  The result holds trivially for $j=1$.
For the base case $j=2$, it follows from Lemma \ref{bigvanishing} that
$$0 = (\Sigma_{1,2\in T}\delta_T)^d \cdot \delta_N^{d(n-2)-1} = \delta_{2}^d \cdot \delta_N^{d(n-2)-1} + \delta_N^d\delta_N^{d(n-2)-1}.$$
Suppose that the result holds for integers less than $j$, and consider
$$0=(\Sigma_{1,j\in T}\delta_T)^d \cdot \delta_2^d \cdots \delta_{j-1}^d \cdot \delta_N^{d(n-j)-1}.$$
If $1,j\in T$, then $\delta_T\delta_{j-1}=0$ unless $\{1,\ldots,j\}\subset T$
 by Lemma \ref{Chowlemma}.  
Nonzero terms in the expansion can only involve 
 $\delta_T= \delta_j,\delta_N$ by Lemma \ref{bigvanishing}; moreover, terms
involving both $\delta_j$ and $\delta_N$ vanish.  Therefore 
$$
\delta_2^d \cdots \delta_{j}^d \cdot \delta_N^{d(n-j)-1}=-\delta_2^d \cdots
\delta_{j-1}^d \cdot \delta_N^{d(n-(j-1))-1}=-(-1)^{j-2}
\delta_N^{n-2}$$
as needed, where the final equality holds by the inductive hypothesis.
\end{proof}

\begin{proposition}\label{prop1.1} For  $1\leq j\leq n$, if
$1\leq i\leq j$ and $1\le k \le j-i$, then
$$\delta_2^d \cdots \delta_{i}^d \cdot \delta_{i+k}^{dk} \cdot
\delta_{i+k+1}^d \cdots \delta_{j}^d \cdot \delta_N^{d(n-j)-1} =
(-1)^{j-k}\delta_N^{d(n-1)-1}.$$
\end{proposition}
\begin{proof}
\noindent
We prove the result by induction on $j-i$. The base case $j-i=1$ so
that $k=1$ is Lemma \ref{squish}.  Suppose that
the statement holds for integers less than $j-i$.
Consider the identity:
\[
0=(\Sigma_{1,i+1\in T}\delta_T)^d \cdot \delta_i \cdot \delta_{i+k}^{d(k-1)} \cdot
\delta_{i+k+1}^d \cdots \delta_j^d\cdot\delta_N^{d(n-j)-1}
\]
Lemma \ref{Chowlemma} implies that  
nonzero contributions involve only $\delta_T$ with
$\{1,\ldots,i+1\}\subset T$.  Terms involving $\delta_T$ with $|T|>j$ 
vanish by Lemma \ref{bigvanishing}, those involving $\delta_T$ with 
$i+1<T|<i+k$ vanish by Lemma \ref{prop1.2}, as do those involving
$\delta_{i+k+1},\ldots, \delta_j$.  Terms involving
both $\delta_{i+1}$ and $\delta_{i+k}$
vanish by Lemma \ref{prop1.2}.
Hence
$$0= \delta_i \cdot \delta_{i+1}^d \cdot
\delta_{(i+1)+(k-1)}^{d(k-1)} \cdot \delta_{i+k+1}^d \cdots
\delta_j^{d(n-j)-1}+\delta_i \cdot \delta_{i+k}^{d(k)} 
\cdot \delta_{i+k+1}^d \cdots \delta_j^d\delta_N^{d(n-j)-1}.$$
Multiplying by $\delta_2^d \cdots \delta_{i-1}^d\delta_i^{d-1}$
and applying the inductive hypothesis gives the result:
\begin{multline*}
\delta_2^d \cdots \delta_i^{d} \cdot \delta_{i+k}^{d(k)} \cdot
\delta_{i+k+1}^d \cdots \delta_j^d\cdot\delta_N^{d(n-j)-1} \\
= -\delta_2^d \cdots \delta_{i+1}^{d} \cdot
\delta_{(i+1)+(k-1)}^{d(k-1)} \cdot \delta_{i+k+1}^d \cdots
\delta_j^d\cdot\delta_N^{d(n-j)-1}
=-(-1)^{j-(k-1)}\delta_N^{d(n-1)-1}.
\end{multline*}
\end{proof}

\begin{proposition} \label{top_intersection}
$\displaystyle \delta_{N}^{d(n-1)-1} = (-1)^{d(n-1) - 1}.$
\end{proposition}

Before proving this proposition, we give a proof of Theorem \ref{pairing}.

\begin{proof} [Proof of Theorem \ref{pairing}]
Let $T \subsetneq N$ with $|T|\ge 2$.  Without loss
of generality, say $T = \{1, \ldots, j\}$.  Then $\delta_T \cdot C_T$ is
equal to
\[
 \delta_{j}^{d(j-1)} \cdot \delta_{N}^{d(n-j)-1}
=-\delta_N^{d(n-1)-1} =  (-1)^{d(n-1)}
\]
by Proposition
\ref{prop1.1} for $i=1$, $k=j-1$ and Proposition \ref{top_intersection}.

Let $S\subsetneq N$.  If $\emptyset\neq S\cap T\subsetneq S,T$, 
then $\delta_S\cdot C_T =0$ by Theorem \ref{Chowring}. 
If $S\cap T=\emptyset$, $s\in S$, then 
$$0=(\Sigma_{1,s\in T'}\delta_{T'})^d \cdot \delta_S \cdot \delta_j \cdot
\delta_N^{d(n-j-1)-1}=\delta_S \cdot \delta_j \cdot
\delta_N^{d(n-j)-1}$$
since 
nonzero summands only involve $\delta_{T'}$ with $T\cup\{s\}\subset T'$
by Lemma \ref{Chowlemma},
and all these contribute zero unless $T'=N$ by Lemma \ref{bigvanishing}.  
Therefore $\delta_S \cdot C_T =0$.
If $T\subsetneq S$, then $\delta_S\cdot C_T =0$ by Lemma \ref{bigvanishing}. 

Suppose $S \subsetneq T$.  We may assume that 
$S=\{1, \ldots, j-k\}$ with $k\geq 1$. Then
$$\delta_S \cdot C_T = \delta_{j-k} \cdot \delta_{j}^{d(j-1)-1}
\cdot \delta_N^{d(n-j)-1}.$$
It follows from Lemma \ref{prop1.2} that this is zero
when $d(j-1)-1\ge dk+1$, or equivalently when $d(j-k-1)\ge 2$.
We conclude that
$\delta_S \cdot C_T = 0$  if $d =1$ and $|S| \ge 3$, or if $d \ge 2$.
If $d=1$ and $|S|=2$, then 
\[
\delta_S \cdot C_T=\delta_2\delta_j^{j-2}\delta_N^{(n-j)-1} 
= \delta_N^{n-2}=(-1)^{n-2}.
\]
by Proposition \ref{prop1.1} for $i=2,k=j-2$ and
Proposition \ref{top_intersection}.
\end{proof}

It  remains to prove Proposition 
\ref{top_intersection}. Since the proof involves
spaces $T_{d,n}$ for varying $n$, for the remainder
of this section, we use the more precise notation $\delta_{S,n}=[T_{d,n}(S)]$.
In this language, we need to show 
$$\int_{T_{d,n}} \delta_{N,n}^{d(n-1)-1} = (-1)^{d(n-1)-1}.$$
We first establish the following lemmas.

\begin{lemma} \label{base_case}
$$\int_{T_{1,3}} [T_{1,3}(\{1,2,3\})] = -1 {\mbox{ and }}
\displaystyle \int_{T_{d,2}} [T_{d,2}(\{1,2\})]^{d-1} = (-1)^{d-1}.$$
\end{lemma}

\begin{proof} \noindent
From Lemma \ref{base_tdn}, we know that $T_{d,2} \cong \bP^{d-1}$.  As
discussed in the beginning of the section, $[T_{d,2}(\{1, 2\})]$
corresponds to $i^*[D(\{1, 2\})]$, where $i$ is the inclusion $D(\{1,
2\})\hra A^d[2]$. By \cite{Ful:IT}, we have
$$i^*[D(\{1, 2\})] = c_1(N_{D(\{1, 2\})} \A^d[2]) = \cO_{D(\{1,2\})}(-1).$$
Hence $[T_{d,2}(\{1,2\})]=c_1(\cO_{\bP^{d-1}}(-1))$, and the result follows.

Write $T_{1, 3} \cong \bP^1$. As above, for $i$ the inclusion 
$D(\{1, 2, 3\})\hra \A^1[3]$, $i^*[D(\{1, 2, 3\})]$
corresponds to $\cO_{\bP^1}(-1)$ on $T_{1,3}$, which gives the result.
\end{proof}

Let $\pi_{n+1} : T_{d,n+1} \ra T_{d,n}$ be the map which drops the
$(n+1)$st marking as described in Section \ref{construction}.

\begin{lemma} \label{pullback_identities}
$\displaystyle \pi_{n+1}^*(\delta_{N,n}) = \delta_{N,n+1} + \delta_{N^+,n+1}.$
\end{lemma}
\begin{proof}
As in the beginning of the section, let $i : D_{\A^d[n]}(N)
\hookrightarrow \A^d[n]$ be the inclusion, and $p : D_{\A^d[n]}(N) \ra
T_{d,n}$ the vector bundle morphism. Let $\pi = \pi_{n+1}$ and
consider the commutative diagram given by dropping the $n+1$'st point:
\bd
\no{T_{d,n+1}} \arr{s,r}{\pi} 
  \no{D_{\A^d[n+1]}(N^+)} \arr{w,t}{p_+} \arr{s,r}{\pi'} \arr{e,t}{i_+} 
  \no{X[n+1]} \arr{s,r}{\pi''}\\
\no{T_{d,n}} 
  \no{D_{\A^d[n]}(N)} \arr{w,t}{p} \arr{e,t}{i} 
  \no{X[n],}
\ed
\noindent
where we have used the subscripts to distinguish which space the
divisor sits on. By \cite{FM}, Proposition 3.4,
$(\pi'')^*[D_{\A^d[n]}(N)] = [D_{\A^d[n+1]}(N^+)] +
[D_{\A^d[n+1]}(N)]$. Therefore, by commutativity of the diagram, we
have :
\begin{multline*}
\pi^*(\delta_{N,n}) = \pi^*(p^*)^{-1} i^! [D_{\A^d[n]}(N)] =
({p_+}^*)^{-1}{i_+}^*(\pi'')^* [D_{\A^d[n]}(N)] \\ = 
({p_+}^*)^{-1}{i_+}^* \big([D_{\A^d[n+1]}(N^+)] +
[D_{\A^d[n+1]}(N)]\big) = \delta_{N,n+1} + \delta_{N^+,n+1}.
\end{multline*}
\end{proof}

\begin{lemma} \label{induction_step}
Let $\pi = \pi_{n+1} : T_{d,n+1} \ra T_{d,n}$ be the map which drops the
$(n+1)$st marking as described in Section \ref{construction}.  Then
$$\int_{T_{d,n+1}} \delta_{N^+,n+1}^{dn-1} = (-1)^d \int_{T_{d,n}}
\pi_*\Big(\delta_{N,n+1}^d\Big) \cdot \delta_{N,n}^{d(n-1)-1}.$$
\end{lemma}

\begin{proof}
We first note that, solving for $\delta_{N^+,n+1}$ 
in Lemma \ref{pullback_identities} gives
\begin{multline*}
\pi_*(\delta_{N^+,n+1}^d)=(\pi^*\delta_{N,n} - \delta_{N,n+1})^d \\
=\sum_{i=0}^{d}{\textstyle \binom{d}{i}}(-1)^i\pi_*(\delta_{N,n}^i) \cdot
\delta_{N,n+1}^{d-i} 
= (-1)^d \pi_*(\delta_{N,n+1}^d),
\end{multline*}
since $\dim T_{d,n}=d(n-1)$, so that 
$\pi_*(\delta_{N,n+1}^i)=0$ for $i < d$.

Again solving for $\delta_{N^+,n+1}$ 
in Lemma \ref{pullback_identities}, we have
\[
\delta_{N^+,n+1}^{dn-1}=\delta_{N^+,n+1}^d
\cdot \Big(\pi^*\delta_{N,n}-\delta_{N,n+1}\Big)^{d(n-1)-1}.
\]
Since $\delta_{N^+,n+1}^d \cdot \delta_{N,n+1} =0$ by Lemma \ref{bigvanishing},
the summands from the binomial expansion 
vanish for any positive power of $\delta_{N,n+1}$.  Hence
\[
\delta_{N^+,n+1}^{dn-1}=
\delta_{N^+,n+1}^d \cdot (\pi^*\delta_{N,n})^{d(n-1)-1}
= \delta_{N^+,n+1}^d \cdot \pi^*(\delta_{N,n}^{d(n-1)-1}).
\]
The projection formula and the fact above gives the result.
\end{proof}

For the next two lemmas, we follow the notation of Section
\ref{construction}.
\begin{lemma} \label{lem1}
Writing $F_{d,n}^1 = \bP_{T_{d,n}}(\cL_N \oplus V_{T_{d,n}})$, we may
identify $F_{d,n}^1(N) \subset F_{d,n}^1$ with the subbundle
$\bP_{T_{d,n}}(V_{T_{d,n}})$.
\end{lemma}
\begin{proof}
We recall that the functor defining $F_{d,n}^1$ takes a variety $H$ to the
collection of compatible screens $(i_S : \cL_S \ra (V_H)^S/V_H)$ where either
$S \subset N$ or $S = N^+$. Equivalently, examining the proof of Theorem
\ref{fm_blowup}, the data of $i_{N^+}$ amounts to choosing a vector bundle
inclusion $j : \cL_{N^+} \ra \cL_N \oplus V_H$. The morphism $i_{N^+}$ is then
obtained by using $i_N$ to map $\cL_N$ into $(V_H)^N/V$ and identifying
$(V_H)^{N^+}/V_H \cong (V_H)^N/V_H \oplus V_H$.

The subfunctor represented by $F_{d,n}^1(N)$ is defined by requiring
that all the compatibility morphisms $\cL_{N^+} \ra \cL_S$ are zero
for $S \subset N$. By compatibility, it suffices to know
that $\cL_{N^+} \ra \cL_N$ is zero, or that the
morphism $j$ maps $\cL_{N^+}$ entirely inside of $V_H \subset
\cL_{N^+} \oplus V_H$. But this just says that $F_{d,n}^1(N) =
\bP_{T_{d,n}}(V_{T_{d,n}}) \subset \bP_{T_{d,n}}(\cL_N \oplus
V_{T_{d,n}})$, as desired.
\end{proof}

\begin{lemma} \label{lem2}
Let $\rho : T_{d,n+1} \ra F_{d,n}^1$ be the natural projection. Then
$\rho^*([F_{d,n}^1(N)]) = [T_{d,n+1}(N)]$.
\end{lemma}
\begin{proof}
We first note that scheme-theoretically
$\rho^{-1}(F_{d,n}^1(N))=T_{d,n+1}(N)$. To see this, we
note that a morphism from a scheme $H$ to $\rho^{-1}(F_{d,n}^1(N))$ is
given by specifying a collection of compatible screens $(\cL_S \ra
(V_H)^S/V_H)$ for $S \subset N^+$ such that the compatibility morphism
$\cL_{N^+} \ra \cL_N$ is zero. But it is easy to check that this is
precisely equivalent to the morphism being in $T_{d,n+1}(N)$.

Now the pullback $\rho^*([F_{d,n}^1(N)])$ is represented by a class on
the scheme theoretic inverse image. Since $\rho$ is of relative dimension
zero, the pullback is represented by a
multiple of the fundamental class of the inverse image
$[T_{d,n+1}(N)]$. But since $\rho_* \rho^* = id$, this multiple must be
$1$.
\end{proof}

\begin{lemma} \label{small_top_intersection}
With the notation of Lemma \ref{induction_step}, 
$$\pi_*\Big(\delta_{N,n+1}^d\Big) = [T_{d,n}].$$
\end{lemma}

\begin{proof}[Proof of Lemma \ref{small_top_intersection}]
Consider the commutative diagram of natural morphisms
\bd
\no[2]{T_{d,n+1}} \arr{se,t}{\pi} \arr{s,l}{\beta} \\
\no[2]{F_{d,n}^1} \arr{e,t}{\rho}
  \no{F_{d,n}^0 = T_{d,n}.}
\ed
\noindent
By Lemma \ref{lem2}, $[T_{d,n+1}(N)] = \beta^*[F_{d,n}^1(N)]$, and so
\begin{multline*}
\pi_*([T_{d,n+1}(N)]^d) = \pi_* \beta^*([F_{d,n}^1(N)]^d) = \rho_* \beta_*
\beta^* ([F_{d,n}^1]^d) \\ = \rho_* \beta_* \beta^* ([F_{d,n}^1(N)]^d)
= \rho_*([F_{d,n}^1(N)]^d)
\end{multline*}
\noindent
By Lemmas \ref{lem1} and \ref{applem}, we have
$\rho_*([F_{d,n}^1(N)]^d) = [T_{d,n}]$, completing the proof.
\end{proof}

We now use these lemmas to prove the main result of this section,
Proposition \ref{top_intersection}.

\begin{proof}[Proof of Proposition \ref{top_intersection}]
We proceed by induction on $n$. The base
cases $d = 1, n = 3$, and $d > 1,n=2$ are proved in Lemma \ref{base_case}. 
Then
$$\int_{T_{d,n+1}} \delta_{N^+,n+1}^{dn-1} = (-1)^d \int_{T_{d,n}}
\delta_{N,n}^{d(n-1)-1}=(-1)^d(-1)^{dn - 1}=(-1)^{d(n+1)-1},$$
as needed, where the first equality follows from
Lemma \ref{induction_step} and Lemma \ref{small_top_intersection},
and the second equality follows from the inductive hypothesis.
\end{proof}

\subsection{Conjectural Pairing of cycles} \label{conjectural_pairing}

Let $\cS$ be a collection of non-overlapping subsets of $N$. For subsets $S, T
\in \cS$, we use the notation $S \prec T$ to mean that $S \subset T$ and for
every $U \in \cS$ such that $U \subset T$, we have $U \subset S$. Previously we
have been calling this relation ``$S$ is a child of $T$.''

\begin{defn}
We define the following symbols:
\begin{enumerate}
\item
For $S \in \cS$, $ch(S) = \{T | T \prec S\}$.
\item
For $S \in \cS$, $\chi(S) = |S| - \sum_{T \in ch(S)} |T| + |ch(S)| - 1$.
\end{enumerate}
\end{defn}

The conjectural formula is as follows:
$$\left(\prod_{N \neq S \in \cS} \delta_S^{d\chi(S)}\right)\delta_N^{d\chi(N) - 1} =
\pm \delta_N^{dn - d - 1}.$$
This gives the following conjectural pairing: For the cycle $\prod_{S \in \cS}
\delta_S^{n_S}$, where each $n_i > 0$, we conjecture that

$$\left(\prod_{S \in \cS} \delta_S^{n_S}\right) \left(\prod_{N \neq S \in \cS}
\delta_S^{d\chi(S) - n_S}\right)\delta_N^{d\chi(N) - n_N - 1} = \pm \delta_N^{dn - d -
1}.$$
In other words, $\prod_{S \in \cS} \delta_S^{n_S}$ is ``dual'' to the cycle
$$\left(\prod_{N \neq S \in \cS} \delta_S^{d\chi(S) -
n_S}\right)\delta_N^{d\chi(N) - n_N - 1}.$$

\begin{example}
For  $\cS = \{S\}$, the codimension $1$-cycle
$\delta_S$ is dual to $\delta_S^{d(|S| - 0 + 0 - 1)}\delta_N^{d(n - s + 1 - 1) - 1}$.
\end{example}

\begin{example}
For $\cS = \{S, T\}$ with $S \subsetneq T$, the codimension $2$-cycle
$\delta_S \delta_T$ is dual to $\delta_S^{d(|S| - 0 + 0 - 1)}\delta_T^{d(|T| - |S| + 1 -
  1)}\delta_N^{d(|N| - |T| + 1 - 1) - 1}$.
\end{example}

\begin{example}
For $\cS = \{S, T\}$ with $S\cap T =\emptyset$,
$\delta_S \delta_T$ is
dual to $\delta_S^{d(|S| - 0 + 0 - 1)}\delta_T^{d(|T| - 0 + 0 - 1)}\delta_N^{d(|N| - |S| -
|T| + 2 - 1) - 1}$.
\end{example}

\section{Appendix}

\subsection{An intersection-theoretic Lemma}
\begin{lemma} \label{applem}
Let $X$ be a variety, and suppose $\pi : E \ra X$ is a vector bundle
of rank $d+1$ and $V \subset E$ is a rank $d$ subvector
bundle. Consider the inclusion $i : P(V) \hra P(E)$. Then
$\pi_*((i_*[P(V)])^d) = [X]$ in the Chow group of $X$.
\end{lemma}
In the proof of this fact we omit the $i_*$ for notational
convenience.
\begin{proof}
Let us first consider the case where $E = L \oplus V$ for some trivial
line bundle $L$. By \cite{Ful:IT}, Lemma 3.3, $c_1(\cO_E(1)) \cap
[P(E)] = [P(V)]$. Using \cite{Ful:IT}, Proposition 3.1, we have $[X] =
s_0(E) \cap [X]$. But by the definition of the Segre class, we have:
\begin{align*}
[X] = s_0(E) \cap [X] &= \pi_*(c_1(\cO_E(1))^d \cap [P(E)]) \\
&= \pi_*([P(V)]^d \cap [P(E)]) = \pi_*([P(V)]^d)
\end{align*}

\noindent
Now let us consider the general case. By counting dimensions, we know
that $\pi_*([P(V)]^d) \in A_{dim(X)}(X) = \Z [X]$, and so
$\pi_*([P(V)]^d) = a [X]$ for some $a \in \Z$. We show that $a =
1$. Choose an open subvariety $j : U \hra X$ such that $E|_U = L|_U
\oplus V|_U$ for some trivial line bundle $L$. We have a commutative
diagram:
\bd
\no[2]{P(V|_U)} \arr{e,t}{i_U} \arr{s,l}{j_V}
  \no{P(E|_U)} \arr{e,t}{\pi_U} \arr{s,l}{j_E}
    \no{U} \arr{s,l}{j} \\
\no[2]{P(V)} \arr{e,b}{i} 
  \no{P(E)} \arr{e,b}{\pi}
    \no{X}
\ed

\noindent
Now we observe that $j^* \pi_* ([P(V)]^d) = j^*(a[X]) = a[U]$. But by
\cite{Ful:IT} Propositions 1.7 and 2.3(d), this gives:
\begin{multline*}
a[U] = j^* \pi_* ([P(V)]^d) = (\pi_U)_* j_E^*([P(V)]^d) \\ =
(\pi_{U})_* \big((j_E^*[P(V)])^d\big) = (\pi_U)_*([P(V|_U)]^d)
\end{multline*}
But by the first case, we have $(\pi_U)_*([P(V|_U)]^d) = [U]$, which
implies that $a = 1$.
\end{proof}

\subsection{Inverse image of a height $1$ prime ideal sheaf}
\begin{lemma} \label{ht1}
Suppose $f : X \ra Y$ is a dominant morphism of varieties and $\cI
\subset \cO_Y$ is the ideal sheaf of a codimension $1$ subvariety $Z
\subset Y$. If $Y$ is locally factorial then the canonical map of
coherent sheaves $f^*(\cI) \ra f^{-1}(\cI)$ is an isomorphism.
\end{lemma}
\begin{proof}
The result is local, and so we may assume that $X = Spec(S), Y =
Spec(R)$, and $\cI = \til{I}$ for an ideal $I \subset R$. By the
hypotheses we may also assume that $R$ and $S$ are domains, $\phi$ is
injective, and that $R$ is factorial, and therefore by \cite{Eis},
Corollary 10.6, $I = aR$ for some $a \in R$. Let $\phi : R \ra S$ be
the map induced by $f$. Then we need to show that $I \otimes_R S \cong
\phi(I)S$, where the isomorphism is induced by the multiplication
map. Since $I = aR$, this amounts to showing that the morphism $S \ra
S$ induced by multiplication by $\phi(a)$ is injective. But this is
true since $\phi$ is injective and $S$ is a domain.
\end{proof}

\subsection{HI Spaces}

In this section we show that if the ground field is the complex
numbers, then the Chow groups of $T_{d,n}$ are isomorphic to the
homology (and cohomology) groups.

\begin{defn}
A complex algebraic variety $X$ is an \emph{HI space} if the canonical map
$cl_X : A_i(X) \ra H_{2i}(X)$ is an isomorphism.
\end{defn}

\begin{defn} \label{cell_def}
A variety $X$ is {\emph{cellular}} if it is has a filtration $X_0
\subsetneq X_1 \subsetneq \cdots \subsetneq X_n = X$, such that $X_i \setminus
X_{i-1}$ is a disjoint union of affine spaces $\A^i$.
\end{defn}

If $X$ has a cellular decomposition then it is an HI space. This fact
can be found in \cite{Ful:IT}, example 19.1.11(b).  Also, if $Y
\subset X$ is a closed subscheme, $U = X \setminus Y$ the open
complement, and $cl_Y$, and $cl_U$ are both isomorphisms, then so is
$cl_X$. That is, $Y, U$ both HI implies $X$ HI. This is in
\cite{Ful:IT}, example 19.1.11(a). It is also easy to check that if
$X$ is a projective bundle over $Y$, and $Y$ is an HI space, then so
is $X$.

We recall the following fact:
\begin{theorem}[\cite{Keel}]
Suppose $Y$ is a closed subvariety of $X$, and $X$ and $Y$ are both
HI. Then The blowup of $X$ along $Y$ is also HI.
\end{theorem}

From this, our description of the space $T_{d,n}$ as a blowup of a projective
bundle shows that:
\begin{corollary}\label{tdn_HI}
$T_{d,n}$ is an HI space.
\end{corollary}

We similarly have:
\begin{corollary}
Suppose that if $X$ is an HI space. Then so is the Fulton-MacPherson
configuration space $X[n]$.
\end{corollary}

\bibliographystyle{alpha}
\bibliography{citations}

\end{document}